\documentclass[final,3p]{elsarticle}
\biboptions{sort&compress}

\usepackage[bitstream-charter]{mathdesign}
\usepackage[T1]{fontenc}

\usepackage[english]{babel}
\usepackage[colorlinks=true]{hyperref}

\usepackage{caption}
\usepackage{subcaption}

\usepackage{amsmath,amsthm}
\usepackage{mathtools}
\usepackage{commath}
\usepackage{enumerate}

\usepackage{graphicx}
\usepackage{booktabs}
\usepackage[nameinlink]{cleveref}

\usepackage{bm}
\usepackage{stmaryrd}
\SetSymbolFont{stmry}{bold}{U}{stmry}{m}{n}
\DeclareMathAlphabet{\pazocal}{OMS}{zplm}{m}{n}
\DeclareMathAlphabet{\pazocalbf}{OMS}{cmsy}{b}{n}
\usepackage{etoolbox}

\preto\subequations{\ifhmode\unskip\fi}
\theoremstyle{plain} \newtheorem{theorem}{Theorem}

\theoremstyle{definition} \newtheorem{definition}[theorem]{Definition}

\theoremstyle{remark} 
\newtheorem{remark}{Remark}

\def\eps{\varepsilon}
\let\phivar\phi
\def\phi{\varphi}

\newcommand{\hcalA}{\hat{\cal A}}

\newcommand{\RR}{\mathbb{R}}

\newcommand{\xb}{\bm{x}}
\newcommand{\xbref}{\hat{\xb}}
\newcommand{\nb}{\bm{n}}

\newcommand{\Cr}{\pazocal{C}}
\newcommand{\CR}{\partial\Cr}

\renewcommand{\O}{\Omega}
\newcommand{\Oref}{\widehat{\O}}

\newcommand{\FL}{\O_f}
\newcommand{\FLref}{\Oref_f}
\newcommand{\SO}{\O_s}
\newcommand{\SOref}{\Oref_s}
\newcommand{\IN}{\Gamma}
\newcommand{\INref}{\widehat{\Gamma}}

\newcommand{\phiref}{\hat{\bm{\phivar}}}
\newcommand{\xiref}{\hat{\xi}}
\newcommand{\zetaref}{\hat{\zeta}}
\newcommand{\psiref}{\hat{\bm{\psi}}}

\newcommand{\fb}{\bm{f}}
\newcommand{\fbref}{\hat{\fb}}

\newcommand{\ub}{\bm{u}}
\newcommand{\ubref}{\hat{\ub}}

\newcommand{\vb}{\bm{v}}
\newcommand{\vbref}{\hat{\vb}}

\newcommand{\pref}{\hat{p}}
\newcommand{\wb}{\bm{w}}

\newcommand{\sigb}{\bm{\sigma}}
\newcommand{\sigbref}{\hat{\sigb}}
\newcommand{\eb}{\bm{e}}

\newcommand{\Vb}{\bm{V}}
\newcommand{\Vbref}{\hat{\Vb}}
\newcommand{\Vref}{\hat{V}_{\theta}}

\newcommand{\hmax}{h_\text{max}}
\newcommand{\hcr}{h_\text{cr}}

\DeclareMathOperator{\COD}{cod}

\DeclareMathOperator{\diver}{div}
\DeclareMathOperator{\inter}{int}

\newcommand{\restr}[2]{{\left.\kern-\nulldelimiterspace#1\right|_{#2}}}

\journal{Journal Name}

\begin{document}

\begin{frontmatter}
\title{A thermo-flow-mechanics-fracture model coupling a phase-field interface approach and thermo-fluid-structure interaction}

\author[1]{Sanghyun Lee}
\ead{slee17@fsu.edu}
\author[2]{Henry von Wahl}
\ead{henry.von.wahl@uni-jena.de}
\author[3]{Thomas Wick}
\ead{thomas.wick@ifam.uni-hannover.de}

\affiliation[1]{organization={Department of Mathematics, Florida State University},
            addressline={1017 Academic Way},
            city={Tallahassee},
            postcode={32306-4510},
            state={FL},
            country={USA}}

\affiliation[2]{organization={Friedrich-Schiller-Universität, Fakultät für Mathematik und Informatik},
            addressline={Ernst-Abber-Platz 2},
            city={Jena},
            postcode={07743},
            country={Germany}}

\affiliation[3]{organization={Leibniz Universit\"at Hannover, Institut f\"ur Angewandte Mathematik},
            addressline={Welfengarten 1},
            city={Hannover},
            postcode={30167},
            country={Germany}}

\begin{abstract}
Geothermal energy, a promising renewable source, relies on efficiently utilizing geothermal reservoirs, especially in Enhanced Geothermal Systems (EGS), where fractures in hot rock formations enhance permeability. Understanding fracture behavior, influenced by temperature changes, is crucial for optimizing energy extraction. To address this, we propose a novel high-accuracy phase-field interface model integrating temperature dynamics into a comprehensive hydraulic-mechanical approach, aiming for a thermo-fluid-structure interaction representation. Therein, the key technical development is a four-step algorithm. This consists of computing the fracture width, reconstructing the sharp interface geometry, solving the thermo-fluid-structure interaction (TFSI) problem, and employing a phase-field approach coupled to the temperature and pressure from the TFSI problem. By coupling temperature-hydraulic-mechanical processes with our newly proposed high-accuracy phase-field interface approach, we investigate how temperature impacts fracture width values, which are crucial for permeability in EGS reservoirs. Through this model and three different numerical simulations, we aim to provide an approach to deepen understanding of the complex interplay between temperature, mechanical deformation, and permeability evolution. Therein, we substantiate our formulations and algorithms through mesh convergence results of crack width and total crack volumes for static fractures, and crack lengths in the case of propagating fractures.

\end{abstract}

\begin{keyword}
phase-field \sep fracture \sep thermo-hydro-mechanics \sep thermo-fluid-structure interaction \sep smeared interface \sep sharp interface
\end{keyword}

\end{frontmatter}

\section{Introduction}

Geothermal energy is a promising renewable energy source, offering sustainable power generation with minimal environmental impact. Understanding the behavior of geothermal reservoirs is crucial for efficient and sustainable exploitation of this resource. In addition, Enhanced Geothermal Systems (EGS) represent a frontier in geothermal energy development, where the creation and stimulation of fractures within hot rock formations are essential for enhancing permeability and facilitating fluid circulation~\cite{moore2019utah,olasolo2016enhanced,lu2018global}. In such systems, understanding the behavior of fractures and their response to changes in temperature is essential to optimize reservoir performance and energy extraction efficiency~\cite{mcclure2014investigation,caulk2016experimental}.

The crack opening displacement or fracture width values within these fractured reservoirs serve as critical permeability indicators, directly influencing fluid flow rates and heat transfer capabilities~\cite{zhang2022thermal,donahue1972crack,burdekin1966crack,zhou2022thermal}. However, accurately predicting these fracture width values requires comprehensive models that account for the coupled effects of temperature, hydraulic processes, mechanical deformation, and fracture propagation.

In numerous works over the last two decades, the phase-field fracture method has been shown to be one of the most effective approaches due to the existence of the diffusive zone~\cite{BourFraMar00, KuMue10, MieWelHof10a, MieWelHof10b, BoVeScoHuLa12, AmGeraLoren15, ARRIAGA201833, SARGADO2018458, WheWiLee20}. Furthermore, monographs and extended papers of phase-field methods for fracture propagation with multi-physics applications include~\cite{BourFraMar08, AmGeraLoren15, WNN20, Wi20, HEIDER2021107881, DiLiWiTy22,WheWiLee20,LeeWheWi16,lee2016phase}. These methods are in contrast to sharp interface approaches such as XFEM~\cite{MDB99,zi2003new,duarte2001generalized,FB10}. However, coupling different physical phenomena may be of interest near or across the interface. The challenge with the phase-field approach lies in the accurate modeling of fundamental physical principles, due to the inherently diffusive nature of the fracture zone, which complicates precise crack boundary localization for modeling interface-related physics. As a result, the phase-field method, while powerful for propagation and representing fracture patterns in two and three dimensions, requires careful consideration when applied to problems where complex physics at the interface between the fluid-filled fracture and the surrounding solid are critical to the fracture width~\cite{LWW17,YNK20}. This is particularly the case when interactions between the fluid and the surrounding solid and temperature effects need to be included. This results in a classical fluid-structure interaction problem with thermal effects in our case, where a sharp interface is required. As a result, the phase-field method, while powerful, requires careful consideration when applied to problems where the exact fracture width plays a critical role.

To address this issue, we build upon a split approach introduced in \cite{WaWi24_RINAM}. This approach combines a diffusive interface phase-field model for fracture dynamics and an interface resolving fluid-structure interaction problem. To combine these two somewhat opposing approaches, we consider a geometry reconstruction approach \cite{WaWi23_CMAME}. To reconstruct the geometry of the open fluid-filled fracture, we use the crack opening displacements or the fracture width to describe the interface between the fluid and the solid. This gives a flexible method for switching between the interface-capturing phase-field method and the interface-tracking fluid-structure interaction approach.

In this paper, we extend the algorithm to consider temperature dynamics into a comprehensive hydraulic-mechanical model using a phase-field approach. By coupling the temperature-hydraulic-mechanical (THM) processes with the phase-field approach, we aim to provide a more realistic representation of fluid-structure interaction processes within geothermal reservoirs, particularly in the context of EGS. Other works considering THM related to phase-field and porous media, include~\cite{NoiiWi19,HEIDER2018116,NgHeiMa23,SUH2021114182,Dai2024,LIU2024117165,lee4920845phase}. Specifically, we model the flow-temperature part through a Boussinesq approximation \cite{LoBo96}. The temperature then enters into the solid equation via the stress tensor; see, e.g., \cite{FARHAT1991349}, \cite[Section 11.4.2]{Wi20} and \cite{Coussy2004}. The THM process is modeled and simulated in the reconstructed domain, which fully resolves the fluid reservoir. However, this differs from the domain in which classical phase-field fractures are considered. The latter is a problem usually posed using a smeared zone of a very thin fracture. This contrasts with the resolved fluid-filled reservoir considered for the THM problem. Therefore, to take quantities, such as the fluid pressure and the temperature from the THM problem, and couple these to a phase-field model, we must derive a novel phase-field fracture model. Consequently, this allows us to switch back to the phase-fracture model on the reconstructed geometry and directly use information from the THM model in the phase-field model.

With this new model at hand, we present a detailed investigation into the influence of temperature on fracture width values using our coupled THM-phase-field model. We hypothesize that temperature variations play a significant role in altering the mechanical properties of the reservoir rock, thereby influencing the opening and closure of fractures and, consequently, permeability enhancement. By incorporating temperature effects into our model, we anticipate gaining deeper insights into the thermal behavior of geothermal reservoirs and its implications for production optimization and reservoir management in EGS.

The remainder of this paper is structured as follows. First, we introduce the system of equations that model the coupled non-isothermal fluid-solid system in \Cref{sec.eqns}. Then, we propose our novel algorithm based on a coupled iteration between the thermo-fluid-structure interaction problem and a phase-field fracture approach. In \Cref{sec.phase_field_model}, we then present a detailed derivation of our fluid-filled fracture phase-field model, which incorporates both pressure and temperature effects at the fracture interface. We then describe the details of how we recover the geometric information from the diffusive phase-field model to reconstruct the sharp interface geometry in \Cref{sec.interface_reconstruction}. Furthermore, we present the weak formulation of the thermo-fluid-structure-interaction problem in this section. In \Cref{sec.numerical_tests}, we then compute several numerical examples to validate and demonstrate the capabilities of the proposed algorithm. Finally, we give some concluding remarks in \Cref{sec.conclusions}.

\section{Governing Equations and Overall Concept}
\label{sec.eqns}

\subsection{Modeling Overview}
\label{sec.eqns:subsec.modeling_overview}

Consider the computational domain $\Omega \subset \mathbb{R}^d$ ($d\in \{2,3\}$). 
We assume that this is subdivided into the fluid domain $\Omega_f$ and the solid domain $\Omega_s$, respectively, as shown in \Cref{fig:setup_problems_a}. In addition, the solid domain $\Omega_s$ is considered a porous medium, while the fluid domain $\Omega_f$ is considered a fracture. The concept of our model is to deal with a sharp interface between $\Omega_f$ and $\Omega_s$, while that interface is moved with the help of a phase-field approach.

\subsubsection{Thermo-Fluid-Structure-Interaction problem}
\label{sec.eqns:subsec.thermo_elasticity_flow}
In our domain, we start with a thermo-fluid-structure-interaction (TFSI) problem. For the fluid and its temperature within the fluid domain $\Omega_f$, we employ the Boussinesq equation. This is a widely adopted approximation to address nonisothermal flow phenomena such as natural convection, circumventing the need to solve the complete compressible formulation of the Navier-Stokes equations. This approximation holds true when density fluctuations are minor, thereby diminishing the problem's nonlinearity. It assumes that density fluctuations minimally influence the flow field, except for their contribution to buoyancy forces. 

The displacement is modeled by linear elasticity in the solid domain $\Omega_s$ while accounting for thermal effects. Consequently, the temperature is solved in the fluid and the solid domain $\Omega_s \cup \Omega_f$. In total, we search for the vector-valued velocity $\bm{v}: \Omega_{f} \rightarrow\mathbb{R}^d$, the vector-valued displacement $\bm{u} : \Omega_s \rightarrow \mathbb{R}^d$, the scalar-valued fluid pressure $p : \Omega_f \rightarrow \mathbb{R}$, and the scalar-valued temperature $\theta : \Omega \rightarrow \mathbb{R}$. These are determined through the following set of equations: Find the velocity, pressure, temperature, and displacement $(\bm{v},p,\bm{u},\theta)$, such that
\begin{subequations}\label{eq_Bou}
\begin{align}
  \rho (\bm{v} \cdot \nabla) \bm{v} -
  \nabla\cdot  \sigb_{Bou}
  - \alpha_\theta (\theta-\theta_0) \fb_f &= 0  &&\text{in } \Omega_f, \label{eq_Bou_1}\\
  \nabla\cdot \bm{v} &= 0 &&\text{in } \Omega_f ,\label{eq_Bou_2}\\
  -\nabla \cdot \sigb_R(\bm{u},p,\theta) &= \fb_s &&\text{in } \Omega_s, \label{eq_Bou_3}\\
  -\nabla\cdot  (\bm{\kappa} \nabla \theta)  + \bm{v}\cdot \nabla \theta &= f_{\theta} &&\text{in } \Omega_f\cup \Omega_s.\label{eq_Bou_4}
\end{align}
\end{subequations}

Equations \eqref{eq_Bou_1},\eqref{eq_Bou_2}, and \eqref{eq_Bou_4} correspond to a Boussinesq approximation~\cite{LoBo96,mayeli2021buoyancy} with the stress tensor
\begin{equation}
  \sigb_{Bou} \coloneqq \rho\nu (\nabla\bm{v} + \nabla\bm{v}^T)  - (p-p_0)I,
\end{equation}
where the fluid density is given as $\rho >0 $ and  the fluid viscosity is $\nu>0$. Here, $\alpha_{\theta}$ is the thermal expansion coefficient, and the external force per unit of mass is given as $\fb_f$.

To consider non-isothermal effects in porous media in the solid domain $\Omega_s$, we assume thermo-poroelasticity~\cite{Coussy2004} with the effective solid stress tensor $\sigb_R$ defined as
\begin{equation}\label{eqn.solid.cauchy-stress}
  \sigb_R \coloneqq  2 \mu \eb(\ub) + \lambda tr(\eb(\ub)) I
            - \alpha_B (p-p_0)I
            - 3\alpha_{\theta} K_{dr} (\theta- \theta_0)I,
\end{equation}
with the Lam\'e parameters $\mu,\lambda>0$, the identity matrix
$I\in\mathbb{R}^{d\times d}$ and the linearized strain tensor
\begin{equation*}
  \eb(\ub) \coloneqq  \frac{1}{2} (\nabla \ub + \nabla \ub^T).
\end{equation*}
Moreover, $\alpha_B\in [0,1]$ is Biot's coefficient, $K_{dr} = 2/3\mu + \lambda$ is the bulk modulus, and $\fb_s$ is an external forcing source term. The values $p_0$ and $\theta_0$ are reference values for the pressure and the temperature, for instance obtained at some initial time or a background state. We remark that specifically, $\theta_0$ will play an important role: $\theta - \theta_0 < 0$ means that we inject colder fluid than the existing fluid in the porous media, and $\theta - \theta_0 > 0$ means that we inject warmer fluid. The reference value $p_0$ plays a relatively less important role in this work (since we only consider the injection) and is set to $p_0 = 0$.

Finally, for the temperature, we consider steady state convection-diffusion heat transfer in the entire domain $\Omega$ where $\bm{\kappa}$ is the heat conductivity coefficient, with $\bm{\kappa}_{| \Omega_f} = \bm{\kappa}_f$ and $\bm{\kappa}_{| \Omega_s} =\bm{\kappa}_s$, and $f_\theta$ is an external forcing source term.

\subsubsection{Interface and Boundary Conditions}
The coupling of $\Omega_s$ and $\Omega_f$ involves integrating the solid domain equations and the fluid domain equations on the interface $\CR$, resulting in a thermo-fluid-structure interaction (TFSI) problem. Here, the sharp interface between $\Omega_s$ and $\Omega_f$ is often referred as the fracture interface in our setup and defined as
$$
  \CR \coloneqq  \overline{\Omega}_f \cap \overline{\Omega}_s.
$$
First, we  define the following notations
$$
  p_{s} = \restr{p}{\Omega_s}, \
  p_{f} = \restr{p}{\Omega_s}, \quad\text{and}\quad
  \theta_{s} = \restr{\theta}{\Omega_s}, \
  \theta_{f} = \restr{\theta}{\Omega_f},
$$
to specify the pressure and the temperature values for each subdomains. Furthermore, we propose 
continuity of the temperature and pressure, as well as continuity of the normal stresses interface conditions
\begin{align*}
  p_{s} = p_{f}, \quad
  \theta_{s} = \theta_{f}, \quad
  \bm{\kappa}_s \nabla \theta_s  \nb =  \bm{\kappa}_f \nabla \theta_f \nb, \quad
\sigb_{R}\nb = \sigb_{\FL}\nb,
\end{align*}
where the normal vector $\nb$ points into the fluid domain $\Omega_f$ (fracture region).
We note that the stress tensor $\sigb_{\FL}$ in $\Omega_f$ is defined as
$$
  \sigb_{\FL} \coloneqq  -(p-p_0)I - \alpha_{\theta}(\theta-\theta_0)I,
$$
where we neglect the displacement in the fluid (fracture) domain. The phase-field formulations of the interface conditions are further discussed  in~\Cref{sec.phase_field_model}.

Finally, we assume the domain boundary $\partial \Omega  \coloneqq \partial\Omega_s \setminus \CR$ contained in the solid boundary, and the system is supplemented by the following (outer) boundary conditions
\begin{align*}
&&&& \bm{u} = \bm{u}_D, \ p = p_D, \ \theta &= \theta_D &&\text{on } \partial\Omega,&&\\
  &&&& \bm{\kappa}_s \nabla \theta \cdot  \nb &= 0 &&\text{on } \partial\Omega,
\end{align*}
where $\bm{u}_D, p_D, \theta_D$ are corresponding Dirichlet boundary conditions for the displacement, the pressure, and the temperature, respectively.

\begin{figure}
  \centering
  \begin{subfigure}[b]{0.35\textwidth}
    \centering
    \includegraphics{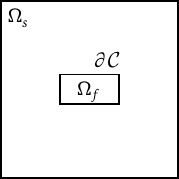}
    \caption{TFSI Domain}
    \label{fig:setup_problems_a}
  \end{subfigure}
  \begin{subfigure}[b]{0.6\textwidth}
    \centering
    \includegraphics{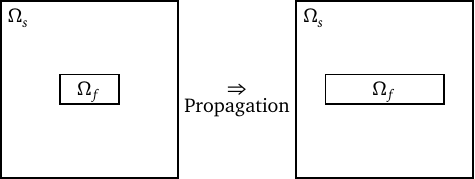}
    \caption{PFF Domain}
    \label{fig:setup_problems_b}
  \end{subfigure}

  \caption{(a) Setup for the thermo-fluid-structure interaction (TFSI) problem  and (b) the change of $\Omega_f$ due to the fracture propagation problem from the phase-field fracture (PFF) approach }
  \label{fig:setup_problems}
\end{figure}

\subsubsection{Phase-Field Fracture Problem}
As we consider the domain $\Omega_f$ as a propagating fracture, we need a method
to compute its changing width and changing length.
In this work, we account for dynamic changes in the subdomains $\Omega_s$ and $\Omega_f$ by considering fracture propagation due to variations in pressure and temperature.
Here, we employ the phase-field fracture (PFF) approach to track the propagation of fractures. Consequently, the fracture (fluid) domain $\Omega_f$ evolves due to the propagation of fractures and the variation in their width and length; see \Cref{fig:setup_problems_b} for an illustration.

The PFF problem not only help to propagate the fracture and tracks the change of $\Omega_f$, but also provides the fracture width (or crack opening displacement) values to create a sharp geometry representation of the fracture. Fracture width, also known as fracture aperture or crack opening displacement (COD), refers to the perpendicular distance between the two opposing faces of a fracture. In the context of geological formations and fluid-structure interactions, it is a critical parameter that influences the flow of fluids through the fracture. The width of a fracture can change over time due to various factors such as stress, pressure/temperature changes, and the propagation of the fracture itself. In modeling and simulations, accurately determining the fracture width is essential for predicting the behavior of fluids within the fractured medium and for understanding the mechanical properties of the fractured solid.
Thus, one of the main goal of this work is to couple the phase-field approach with the thermo-fluid-structure interaction approach to accurately assess the fracture width,  the deformation of the fracture, and related physics across the fracture interface.

\begin{figure}
  \centering
  \includegraphics{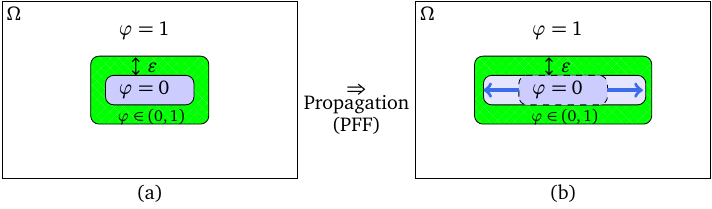}
  \caption{Classical phase-field approach utilizes the diffusive fracture and propagate directly by solving the phase-field problem.}
  \label{fig:phasefield}
\end{figure}

We reconstruct the fractured domain $\Omega_f$ from the phase-field variable $\varphi$ to fully resolve the fracture interface $\CR$ between the fluid and the intact solid.
Here, the scalar-valued phase-field function, $\varphi: \Omega \rightarrow [0,1]$, acts as an indicator function. For example, the fracture domain is defined where $\varphi = 0$, and the intact domain is defined where $\varphi = 1$. The sharp fracture interface $\CR$ becomes a diffusive area/volume domain because the phase-field has a diffusive zone where $\varphi \in (0,1)$ with a characteristic length scale $\varepsilon$; see \Cref{fig:phasefield} (a). The phase-field approach solves the fracture problem to propagate the fracture by tracking the phase field values to simulate the fracture propagation as illustrated in \Cref{fig:phasefield} (b).

\subsubsection{Representation of Fractures}
\begin{figure}[t]
  \centering
  \begin{subfigure}[b]{0.3\textwidth}
    \centering
    \includegraphics{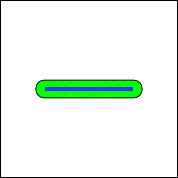}
  \caption{Phase-field fracture (PFF)}
  \label{fig:type_of_fractures_a}
  \end{subfigure}
  \begin{subfigure}[b]{0.3\textwidth}
    \centering
    \includegraphics{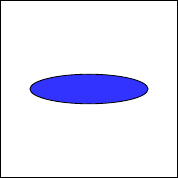}
    \caption{Sharp interface ellipse fracture }
    \label{fig:type_of_fractures_b}
  \end{subfigure}
  \begin{subfigure}[b]{0.3\textwidth}
    \centering
    \includegraphics{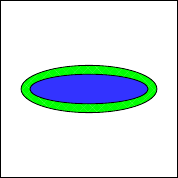}
    \caption{Diffusive interface ellipse fracture }
    \label{fig:type_of_fractures_c}
  \end{subfigure}
  \caption{Illustration of different types of fractures}
  \label{fig:type_of_fractures}
\end{figure}

For the coupling of the TFSI and PFF problems, careful definitions of the fracture are required. The TFSI problems involve a sharp interface $\CR$ between $\Omega_f$ and $\Omega_s$, whereas the PFF problems consider a diffusive interface.
First, the \textit{classical diffusive phase-field fracture (PFF)} is illustrated in \Cref{fig:type_of_fractures_a}. This is the diffusive fracture obtained by solving the classical PFF problem. We note the diffusive zone where $\varphi \in (0,1)$ around the thin fracture zone (where $\varphi = 0$).
Secondly, \Cref{fig:type_of_fractures_b} presents the sharp interface, ellipse-shaped fracture. This fracture is obtained by computing the COD values from the PFF shown in \Cref{fig:type_of_fractures_a}, and it has no diffusive zone. We utilize this fracture to solve the TFSI problem.
Finally, the sharp interface ellipse fracture is converted back to a diffusive interface ellipse fracture by employing the PFF problem.

\subsection{Overall Coupled Algorithm}
\label{sec.eqns:subsec.overall_algorithm}

\begin{figure}
  \centering
  \includegraphics{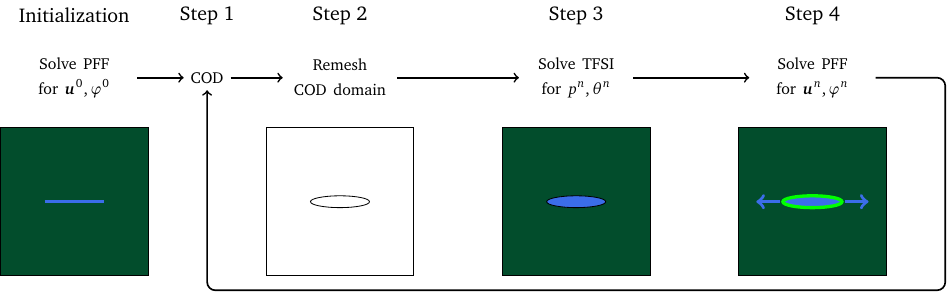}
  \caption{Sketch of the presented algorithm.}
  \label{fig:alg-sketch}
\end{figure}

In this section, we discuss our proposed algorithm, which considers coupling the non-isothermal TFSI problem to the classical phase-field fracture (PFF) problem. The overall concept can be summarized as follows~\cite{WaWi23_CMAME, WaWi24_RINAM}:
\begin{enumerate}\addtocounter{enumi}{-1}
    \item \textbf{Initialization}. This step is only done for the initial time, where we obtain the \textit{classical diffusive phase-field fracture} $\varphi^0$ (initial phase field) by solving the PFF problem with given initial pressure $p^0$ and temperature $\theta^0$.
    \item \textbf{Step 1}. Compute the fracture width (crack opening displacement (COD)) values to create a sharp geometry representation of the \textit{sharp interface ellipse fracture}.
    \item \textbf{Step 2}. Reconstruct the geometry of the open fluid-filled fracture $\Omega_s$ and $\Omega_f$, based on the previously computed COD values.
    \item \textbf{Step 3}. Solve coupled thermo-fluid-structure interaction (TFSI) problem to get $(\bm{v}, p, \theta)$ in $\Omega_f$ and $(\bm{u}, \theta)$ in $\Omega_s$.
    \item \textbf{Step 4}. Given the pressure $p$ and the temperature $\theta$, solve the PFF problem to obtain the displacement and the phase field $(\bm{u}, \varphi)$.  This step is considered to be the prediction of the phase-field fracture, and provides the new fracture domain $\Omega_f$.
    The difference between the initialization and and Step 4 is that the phase-field fracture representation in Step 4 is \textit{diffusive interface ellipse fracture} whereas the initialization is the thin classical PFF.
\end{enumerate}
A sketch of the algorithm can be seen in \Cref{fig:alg-sketch}.
The above algorithm considers multiple couplings between the different variables. First, the temperature and velocity are coupled through the convective and buoyancy terms. Secondly the temperature couples to the solid stress. Finally, the displacement couples to the fluid and temperature equations through the boundary between the fluid and the solid domains.

\section{Governing System for a Phase-Field Fracture Model}
\label{sec.phase_field_model}

In the previous section, we introduced our overall solution concept. Here, we derive the phase-field fracture approach (PFF problem) with the governing equations used in
this work to model the fluid-filled non-isothermal fracture (Initialization and Step 4).
This PFF problem solves the system to obtain the unknown vector-valued displacement field $\ub$ and the scalar-valued phase-field function $\varphi$.

In this section begin with a classical formulation and then present the linearized and regularized formulation used in our implementation. 
Crucially, the latter takes into account that the fluid-pressure driving the fracture is only available inside the open crack from an FSI problem. In the following, starting from quasi-static phase-field fracture modeling based on the original work from~\cite{BourFraMar00}, we recapitulate how pressure interface conditions are included such 
that a pressurized phase-field fracture model is obtained~\cite{MiWheWi13a,MWW19}.

Let us assume we have a fracture $\mathcal{C} \subset \mathbb{R}^{d-1}$ in the domain $\Omega$, then note we have $\Omega_s = \Omega \backslash \mathcal{C}$. To get the correct contributions from traction boundary forces, it is convenient to start from the energy level. Here, traction forces can be described as
\begin{equation}
\label{eqn:energy_functional}
\int_{\partial \Omega_s} \sigb_R(\ub,p,\theta) \nb \cdot \ub \dif s.
\end{equation}
These forces then form the starting point to derive the driving contributions to the phase-field model.

\begin{remark}
For the following, we note that $\mathcal{C} \subset \Omega_f$ in our current setting. However, $\Omega_f \subset \mathbb{R}^d$ will later approximate $\mathcal{C}$ by utilizing the phase-field function. Thus, $\CR = \partial \Omega_f$ will also approximate $\mathcal{C}$, and we obatin $\mathcal{C} \approx \CR$ in the global formulation (\Cref{sec.globalformulation}) with the phase-field function.
\end{remark}

\subsection{Interface Conditions}
\label{subsec.interface.conditions}

The boundary conditions $\partial\Omega$ may be chosen appropriately for a given problem under consideration. For example, we assume homogeneous Dirichlet boundary condition for $\ub$. However, care must be taken to obtain the correct interface conditions is necessary for $\CR$, since $\mathcal{C} \approx \CR$, and $\CR = \overline{\Omega}_s \cap
\overline{\Omega}_f$ becomes the interface between the fluid (fracture) and the
solid (intact) domain. We follow \cite{MiWheWi13a} to model the pressure interface
conditions between the surrounding medium and the fracture, and we refer to \cite{NoiiWi19} and \cite[Section 11.4.2.1]{Wi20} regarding the interface conditions for both the pressure and temperature. The resulting model is a non-isothermal, pressurized phase-field fracture approach. In the following, we provide the mathematical details from prescribing the integral interface conditions and their equivalent formulation as domain integrals.

To include the traction forces \eqref{eqn:energy_functional} in the phase-field model, we assume continuity of the normal stresses on $\CR$. It then follows that~\cite{MiWheWi13a}:
\begin{align}
- \int_{\mathcal{C}} \!\!\!\sigb_{R}\nb {\cdot\ub} \dif s
  &= -\int_{\mathcal{C}} \sigb_{\FL}\nb{\cdot\ub} \dif s \label{eq_int_1} \\
  &= \int_{\mathcal{C}} (p-p_0)\nb\cdot \ub \dif s
   + \int_{\mathcal{C}} \alpha_\theta (\theta-\theta_0) \nb\cdot \ub \dif s \label{eq_int_2} \\
  &=    \int_{\Omega_s}     \!\!\!\! \nabla\cdot ((p-p_0) \ub)\dif\xb 
    -\! \int_{\partial\O} \!\!\!\! (p-p_0)\nb \cdot \ub\dif s 
    +\! \int_{\Omega_s}   \!\!\!\! \alpha_\theta\nabla\cdot ((\theta-\theta_0) \ub)\dif\xb
    -\! \int_{\partial\O} \!\!\!\! \alpha_\theta(\theta-\theta_0)\nb \cdot \ub\dif s \label{eq_int_3}\\
  &= \int_{\Omega_s} (\ub\nabla (p-p_0) + (p-p_0)\nabla\cdot \ub)\dif\xb
    + \int_{\Omega_s} \alpha_\theta (\ub\nabla (\theta-\theta_0)
    + (\theta-\theta_0)\nabla\cdot \ub)\dif\xb,
  \label{eq_int_5}
\end{align}
where Gauss' divergence theorem is applied in \eqref{eq_int_3}, and the
homogeneous Dirichlet conditions $\ub = 0$ on $\partial\Omega$ are employed
in \eqref{eq_int_5}.

\subsection{\texorpdfstring{Global Phase-Field Fracture Formulation}{}}
\label{sec.globalformulation}

To transform integrals from subdomains to the global domain $\Omega$, we follow the standard technique in phase-field fracture, and introduces a degradation function, given by
\begin{equation*}
  g(\varphi)\coloneqq (1-\kappa) {\varphi}^2 + \kappa,
\end{equation*}
with the (small) bulk regularization parameter $\kappa>0$. We note that $g(\varphi)\approx 0$ in the fracture domain (i.e., $\Omega_f$) and $g(\varphi)\approx 1$ in the intact domain (i.e., $\Omega_s$).

Next, phase-field models start from lower-dimensional fractures, described by
the Hausdorff measure of the fracture $\mathcal{C}$. This is then approximated by an Ambrosio-Tortorelli type functional~\cite{AmTo90,AmTo92}:
\begin{equation}
    G_c H^{d-1}(\mathcal{C}) = G_c \int_\Omega \dfrac{1}{2\varepsilon} (1-\varphi)^2 + \dfrac{\varepsilon}{2}(\nabla \varphi)^2 \dif\xb,
\end{equation}
with the critical energy release rate $G_c> 0$. Moreover, the classical phase-field fracture model does not allow an open crack to reseal. Thus, the phase field is subject to the crack irreversibility constraint $\partial_t \varphi \leq 0$. In our phase-field model, this continuous irreversibility constraint is approximated through a difference quotient by $\varphi \leq \varphi^{\text{old}}$. Thus, the phase-field fracture problem is often referred to be in a quasi-static regime. Particularly, the formulation does not contain any time derivatives. Nevertheless, temporal dependence may enter the system through factors such as time-dependent pressure and temperature, and to satisfy the irreversibility constraint.
Due to the quasi-static nature of this problem formulation, we apply some further approximations to arrive at the time-discretized problem.

Let the iteration steps be denoted by $\tau^m$, with index $m\in\mathbb{N}_0$. Then, we have  $\varphi^m \approx \varphi(\tau^m)$. In other works, this iteration is also known as incremental steps, pseudo-time steps, or time steps. We choose our nomenclature to distinguish this from (time) steps $t^n$ with index $n$ used to advance the overall coupled system of phase-field mesh reconstruction and fluid-structure interaction.

Finally, to formulate the  weak form of global phase-field formulation,
we consider the function spaces $W\coloneqq H^1(\O)$,
$\Vb\coloneqq [H^1_0(\O)]^d$ and the convex set
\begin{equation*}
  K\coloneqq \{w\in H^1(\O) |\, w\leq \varphi^{old} \leq 1 \text{ a.e. on }\O\} \cap L^\infty(\O).
\end{equation*}
Then, our proposed non-isothermal, pressurized phase-field fracture problem is given be the following definition.

\begin{definition}\label{form_1}
Let the pressure $p\in W^{1,\infty}(\O)$, the temperature $\theta\in W^{1,\infty}(\O)$,
Dirichlet boundary data $\ub_D$ on $\partial\Omega$, and the initial condition
$\varphi(0)\coloneqq\varphi^0$ be given. Furthermore, let the phase-field
regularization parameter $\eps>0$ and the critical energy release rate $G_c> 0$
be given. We define the interface driven coupled thermo-phase-field fracture
problem as follows. For the iteration steps $m=1,2,3,\dots,M$, find
$(\ub,\varphi)\coloneqq(\ub^{m},\varphi^{m}) \in \{\Vb + \ub_D\} \times K$,
such that
\begin{subequations}\label{eqn.phase-field.form1}
\begin{align}
  \begin{multlined}[b]
  \big(g(\varphi) \; \sigb_R(\ub), \eb( {\wb})\big)_\O
    + \big(g({\varphi}) (p-p_0), \nabla\cdot\wb\big)_\O
    + \big(g({\varphi}) \nabla (p-p_0), \wb\big)_\O\\
    + \big(\alpha_{\theta}g(\varphi) (\theta-\theta_0), \nabla\cdot\wb\big)_\O
    + \big(\alpha_{\theta}g(\varphi)   \nabla(\theta-\theta_0), \wb\big)_\O
  \end{multlined}
  &= 0 ~~\forall \wb\in \Vb,\\
  \begin{multlined}[b]
    (1-\kappa) \big({\varphi} \; \sigb_R(\ub):\eb(\ub), \psi {-\varphi}\big)_\O
    +  2(1-\kappa) \big({\varphi}\;  (p-p_0)\; \nabla\cdot  \ub,\psi{-\varphi}\big)_\O\\
    +  2(1-\kappa) \big({\varphi}\;  \nabla (p-p_0)\; \ub,\psi\big)_\O
    + 2(1-\kappa)\big(\alpha_{\theta}\varphi\; (\theta-\theta_0)\; \nabla\cdot \ub, \psi\big)_\O\\
    + 2(1-\kappa)\big(\alpha_{\theta} \; (\varphi-\theta_0) \; \nabla \theta\, \ub, \psi\big)_\O
    + G_c  \Bigl(\eps \big(\nabla \varphi, \nabla (\psi - {\varphi})\big)_\O
      -\frac{1}{\eps} \big(1-\varphi,\psi{-\varphi}\big)_\O    \Bigr)
  \end{multlined}
  &\geq  0 ~~\forall \psi \in K.
\end{align}
\end{subequations}
\end{definition}
This formulation uses the above interface law formulated as a domain integral
using the Gauss divergence theorem \cite[Section 2]{MiWheWi13a}, as derived in
\Cref{subsec.interface.conditions}.

\subsection{\texorpdfstring{Interface Phase-Field Fracture Formulation}{}}
Now, the above formulation assumes that the temperature and pressure are given in $\Omega=\inter(\overline{\Omega_s \setminus\mathcal{C}})$. However, our aim is to couple the temperature and pressure from a thermo-fluid-structure interaction problem to this phase-field model. Consequently, the pressure will only be available in $\Omega_f$ and the temperature will be defined in $\Omega=\overline{\Omega}_f\cup\overline{\Omega}_s$. Consequently, we need to derive a formulation involving the fracture boundary $\CR=\partial\Omega_f\cap\partial\Omega_s$.

To this end, we recall $\sigb_R$ from \eqref{eqn.solid.cauchy-stress}, and split it into
$$
\sigb_R = \sigb_s - \alpha_B (p-p_0)I - 3\alpha_{\theta} K_{dr} (\theta - \theta_0)I \qquad\text{in } \Omega_s,
$$
where $\sigb_s = 2 \mu \eb(\ub) + \lambda tr(\eb(\ub)) I$ is the linear elasticity part with the displacement.
As in \eqref{eq_int_2}, we do not consider the stress
contributions, such that only the pressure and temperature components interact
from $\Omega_s$ to $\CR$, i.e
\begin{equation}\label{eq_int_6_May_23}
  - \alpha_B (p-p_0)I - 3\alpha_{\theta} K_{dr} (\theta - \theta_0)I \quad\text{in } \Omega_s
\end{equation}
As shown from \eqref{eq_int_1} to \eqref{eq_int_5}, which
transforms the interface integrals to the domain integrals, we perform the similar procedure.

Here, we transform \eqref{eq_int_6_May_23} into an interface integral.
To this end, we work again on the energy level with $\ub$ as variation,
we go backwards the chain and obtain
\begin{equation}\label{eqn.deriv-model:boudary-forces}
- \alpha_B \int_{\CR}  (p-p_0)\nb \ub \dif{s} - 3\alpha_{\theta} K_{dr}\int_{\CR} (\theta- \theta_0)\nb \ub\dif{s}.
\end{equation}
We note that we are now employing $\CR$ instead of $\mathcal{C}$ due to the
given phase-field fracture domain. To transition from the sharp 
interface to the diffusive phase-field representation, we have to include the 
phase-field variable in \eqref{eqn.deriv-model:boudary-forces}. As for the case $\varphi = 0$ the entire integral would vanish,
we add the regularization $\kappa >0$ such that the discrete system matrices remain well-posed.

We then have
\begin{equation}
- \alpha_B \int_{\CR} g(\varphi) (p-p_0)\nb \ub \dif{s} - 3\alpha_{\theta} K_{dr}\int_{\CR} g(\varphi)(\theta- \theta_0)\nb \ub\dif{s}.
\end{equation}
Differentiating in $\ub$ in the direction $\wb$ and in $\varphi$ in the direction $\psi$ yields
\begin{equation*}
  - \alpha_B \int_{\CR} g(\varphi) (p-p_0)\nb \wb \dif{s} - 3\alpha_{\theta} K_{dr}\int_{\CR} g(\varphi)(\theta- \theta_0)\nb \wb\dif{s},
\end{equation*}
and
\begin{equation*}
  - 2(1-\kappa)\alpha_B \int_{\CR} \varphi (p-p_0)\nb \ub \psi\dif{s}
  - 2(1-\kappa) 3\alpha_{\theta} K_{dr}\int_{\CR} \varphi(\theta- \theta_0)\nb \ub \psi\dif{s},
\end{equation*}
respectively. We notice that $\sigb_s$ remains as the domain integral. Thus,
the pressure and temperature contributions from $\sigb_R$ enter in the formulation
as interface integrals on $\CR$ and the solid stress $\sigb_s$ enters as domain
integral contribution. With the above derivation, we finally have the following
non-isothermal, pressurized interface phase-field problem.

\begin{definition}\label{form_2_interface_b}
  Let the data from \Cref{form_1} be given, and $\nb$ denote the unit normal
  vector pointing into the crack. We define the semi-linearized
  interface driven coupled thermo-phase-field problem as follows. For the
  iteration steps $m=1,2,3,\dots,M$, find $(\ub,\varphi)\coloneqq(\ub^m,\varphi^m)
  \in \{\Vb + \ub_D\}\times W$, such that
  \begin{subequations}\label{eqn.form_2_interface_b}
  \begin{align}
    \big(g(\varphi^{m-1})\; \sigb_s(\ub), \eb({\wb})\big)_\O
      \!+\!\int_{\CR}\!\!(1 - \alpha_B)(p - p_0) \nb \cdot \wb \dif{s}
      \!+\!\int_{\CR}\!\!(\alpha_\theta - 3 \alpha_\theta K_{dr})\; (\theta - \theta_0) \nb \cdot \wb \dif{s}
    &= 0
    ~~\forall \wb\in \Vb,\label{eqn.form_2_interface_b.a}\\
    \begin{multlined}[b]
      (1-\kappa) \big({\varphi} \; \sigb_s(\ub):\eb(\ub), \psi\big)_\O
        +  G_c  \Bigl(\eps \big(\nabla\varphi, \nabla\psi\big)_\O -\frac{1}{\eps} \big(1-\varphi,\psi\big)_\O  \Bigr)
        + \big(\gamma(\varphi - \varphi^{m-1})^+,\psi\big)_\O\\
      +  2 (1-\kappa)\int_{\CR}\!(1 - \alpha_B)(p - p_0) \nb \cdot \ub \psi \dif{s}
        +  2 (1-\kappa) \int_{\CR}\!(\alpha_\theta - 3 \alpha_\theta K_{dr})\; (\theta - \theta_0) \nb \cdot \ub \psi \dif{s}
    \end{multlined}
     &=0~~\forall \psi\in W.\label{eqn.form_2_interface_b.b}
  \end{align}
  \end{subequations}
\end{definition}
To the best of our knowledge, a phase-field formulation based on
interface couplings, such as \Cref{form_2_interface_b} has only been used in
\cite{WaWi24_RINAM}. This is because the interface formulation seems to contradict
the phase-field concept, where the interface is not known exactly.
In this paper, we extend the idea in \cite{WaWi24_RINAM} to consider thermal effects.

\begin{remark}[Linearization]
In \Cref{form_2_interface_b} we have relaxed the non-linear behavior in the
first term in \eqref{eqn.phase-field.form1} by using
the approximation $\varphi \approx \varphi^{m-1}$, i.e.,
\begin{equation*}
  g(\varphi)\,\sigb_R(\ub) \quad\mapsto\quad g(\varphi^{m-1})\,\sigb_R(\ub).
\end{equation*}
This follows the extrapolation introduced in~\cite{HWW15} and is numerically
justified for slowly growing fractures~\cite[Chapter 6]{Wi20}.
In the case of fast-growing fractures, this is known to fail~\cite{Wic17} due
to the time lagging errors. In addition, this approximation could introduce the
fix point iteration error which vanishes with the number of iterations.
Alternatively, fully monolithic schemes~\cite{Wic17} or an additional
iteration~\cite{KMW23} must be introduced, to avoid this error.
\end{remark}

\begin{remark}[Penalty method]
The second approximation in \Cref{form_2_interface_b} addresses 
irreversibility constraint. We relax this inequality
constraint by considering a simple penalization, see~\cite{MWT15}
or \cite[Chapter 5]{Wi20}), i.e.,
\begin{equation*}
  \varphi \leq \varphi^{m-1} \quad\mapsto\quad \gamma(\varphi - \varphi^{m-1})^+.
\end{equation*}
Here, $(x)^+ = x$ for $x>0$ and $(x)^+ = 0$ for $x\leq 0$,
and where $\gamma>0$ is a penalty parameter.
\end{remark}

\begin{remark}[Temperature influence on fracture aperture and propagation]
Moreover, we have the following relation for the temperature
interface terms. We have $K_{dr}>1$, and so it holds
\[
(\alpha_\theta - 3 \alpha_\theta K_{dr}) < 0.
\]
This means that the injection of colder fluid than the existing fluid in the
porous media, i.e., $(\theta - \theta_0)<0$, will cause the fracture to
increase in width and length. On the other hand, for warmer water injection,
i.e., $(\theta - \theta_0)>0$, the fracture width will decrease~\cite{TrSeNg13}.
\end{remark}

\begin{remark}[Biot's coefficient $\alpha_B$]
In this work, we assume $\alpha_B = 0$. If the full coupling in poroelasticity
with $\alpha_B =1$ is assumed, one needs to solve for the pressure by utilizing
the poroelastic coupling, then the interface integrals
$\int_{\CR}(1 - \alpha_B)(p - p_0) \nb \cdot \wb \dif{s}$
and $2(1-\kappa) \int_{\CR} (1 - \alpha_B)(p - p_0) \nb \cdot \ub \psi \dif{s}$
vanish. However, in our derivation, we still have the pressure contributions
(from the TFSI problem) at the interface as $p$ still contributes to $\sigb_R$.
\end{remark}

\section{Interface Reconstruction, Remeshing, and Coupled PFF-FSI Framework}
\label{sec.interface_reconstruction}

With the phase-field fracture model presented above, we may compute a fracture 
if the pressure and temperature data are given. However, to obtain these 
quantities from the considered TFSI problem~\eqref{eq_Bou}, we must first 
obtain the geometry in which the problem is posed.
Following our previous work in \cite{WaWi24_RINAM}, we use the geometry
reconstruction approach of the fluid-solid interface presented in \cite{WaWi23_CMAME}.
In this approach, we construct a fitted mesh of the open crack and the
surrounding solid are then able to pose our TFSI
problem on this geometry using an interface tracking approach discussed next.

\subsection{Step 1 and Step 2: Computing COD and Remeshing}
The geometry reconstruction is based crack opening displacement (COD), or aperture
of the crack~\cite{WaWi23_CMAME}. This can be computed by
\begin{equation}\label{eqn:cod-definition}
 \COD(\xb) = \llbracket \ub\cdot\nb \rrbracket
 \simeq \int_{\ell^{\xb,\bm{v}}} \ub(\xb) \cdot \nabla \varphi(\xb) \dif s,
\end{equation}
where $\ell^{\xb,\bm{v}}$ is a line through $\xb$ along the vector $\bm{v}$
\cite{CHUKWUDOZIE2019957}. See also \cite[Proposition~83]{Wi20} for a
simplification of the formula when the crack is aligned with a Cartesian axis.

\begin{figure}
  \centering
  \includegraphics[height=4cm]{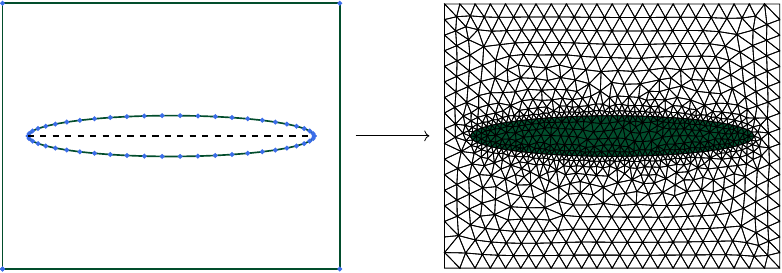}
  \caption{Sketch of re-meshing procedure. The dashed line represents the known
    center line of the crack, the blue points the computed crack opening
    displacements on the crack interface, which are connected by the green line
    segments to form the reconstructed geometry approximation.}
  \label{fig:re-meshing}
\end{figure}

We assume that the centerline of the crack is known, i.e., the line such that
half the fracture width lies on either side of this line, c.f. the dashed line
on the left of \Cref{fig:re-meshing}. With the knowledge of
this line, the crack opening displacements give a set of points on the boundary
of the open crack domain, c.f the blue points on the left of \Cref{fig:re-meshing}.
These can be connected by line segments or higher-order splines, c.f. 
the green line segments on the left of \Cref{fig:re-meshing}. This then forms an
approximation of the crack interface. This geometry can the be remeshed
using an automated meshing tool, resulting in an appropriate mesh for the
a finite element based fluid-structure interaction solver, c.f. the right of 
\Cref{fig:re-meshing}.
We further assume that the number of CODs computed is sufficiently large, such
that the geometry of the open crack is sufficiently well resolved. This can be
achieved by computing $\pazocal{O}(h^{-1})$ crack opening displacements
along the centerline of the crack.

This re-meshing approach has the advantage
that we can consider propagating and merging cracks, while driving the
crack using accurate quantities from a fluid-structure interaction model
with a resolved interface between the solid and the fluid-filled crack.

\begin{remark}
  Computing the COD can become numerically unstable near the tips of the crack, c.f. \Cref{fig.ex2.cod-smoothing} below. This is especially the case when iterating between the phase-field fracture and thermo-fluid-structure interaction problems. To avoid this issue, it can become necessary to preprocess the COD values to smooth out oscillations in the COD. Details of this are given below in \Cref{sec.numerical_approx:subsec.ex2}.
\end{remark}

\subsection{Step 3: Stationary Thermo-Fluid-Structure Interaction in the Reference Configuration}
\label{sec.eqns:subsec.fsi_in_ale}

To obtain a weak formulation suitable for a finite element simulation,
we model the thermo-fluid-structure interaction (TFSI)
problem from \Cref{sec.eqns:subsec.thermo_elasticity_flow}
in arbitrary Lagrangian-Eulerian (ALE) coordinates. This uses
variational monolithic coupling in a reference configuration
\cite{HrTu06a,Du07,Wi11_phd,Ri17_fsi}.
For our model, we assume a stationary flow inside the fluid-filled crack and,
therefore, consider the stationary Navier-Stokes equations for the fluid.
This goes beyond the previously considered Stokes flow in
\cite{WaWi23_CMAME,WaWi24_RINAM}.

For the TFSI problem, let the domain $\O\subset\RR^d$ be divided into
a $d$-dimensional fluid domain $\FL$, a $d$-dimensional solid domain
$\SO$ and a $d-1$-dimensional interface $\IN$ between the two, such that
$\O=\FL\dot\cup\IN\dot\cup\SO$. Furthermore, we also require these domains
in a reference configuration, which we denote by $\hat{\O}, \FLref,
\SOref$ and $\INref$. Similarly, we denote by $\vbref, \pref, \ubref$ and
$\xbref$ the velocity, pressure, deformation and coordinates in the reference
configuration.
In the present setting of a fluid-filled crack, the fluid domain
is the interior of the crack $\FL=\Cr$, the solid is the intact medium
$\SO$ and the interface is the crack boundary $\IN=\CR$.

Formulating the problem in the domains $\SOref$ and $\FLref$ leads to the
well-established formulation in ALE coordinates~\cite{HuLiZi81,DoGiuHa82}. To
obtain a monolithic formulation, we need a transformation $\hcalA_f$ from the
reference configuration to the physical domain in the fluid-domain. This
transformation is given on the interface $\INref$ by the structure displacement:
\begin{equation*}
  \hcalA_f(\xbref,t)\big|_{\INref} = \xbref+\ubref_s(\xbref,t)\big|_{\INref}.
\end{equation*}
On the outer boundary of the fluid domain
$\partial\FLref\setminus\INref$, it holds $\hcalA_f=\text{id}$.
Inside $\FLref$, the only requirement on the transformation is that it should
be as smooth and regular as possible. To this end, we use a harmonic
extension of $\ubref_s|_{\SOref}$ to the fluid domain $\FLref$ and
define $\hcalA_f\coloneqq \text{id}+\ubref$ on $\FLref$. That is,
$id(\xbref) = \xbref$, i.e., $\hcalA_f\coloneqq \text{id}+\ubref
= \hcalA_f(\xbref,t) = \text{id}(\xbref)+\ubref(\xbref,t)$,
such that
\begin{equation*}
  (\hat \nabla \ubref_f,\hat \nabla\psiref)_{\FLref} = 0,\quad
  \ubref_f=\ubref_s\text{ on }\INref,\quad
  \ubref_f=0\text{ on }\partial\FLref\setminus\INref.
\end{equation*}
Consequently, we define a continuous deformation $\ubref$ on all $\Omega$,
which coincides with the solid deformation in $\SOref$ and gives the
appropriate transformation in $\FLref$. Skipping the subscripts and
because $\hcalA_f$ coincides with the solid transformation $\hcalA_s$,
we define on the entire domain $\Oref$:
\begin{equation*}
  \hcalA(\xbref, t)\coloneqq \xbref + \ubref(\xbref, t),\quad
  \hat F(\xbref, t) \coloneqq \hat\nabla \hcalA= I+\hat\nabla \ubref(\xbref, t),\quad
  \hat J\coloneqq \text{det}(\hat F).
\end{equation*}

With this transformation into the reference configuration at hand, we
present weak formulation of the stationary thermo-fluid-structure
interaction problem, see also~\cite{RiWi10} for the formal derivation.

\begin{definition}[Stationary thermo-fluid-structure interaction]
\label{def:fsi:ale:stationary}
Let $\Vbref$ be a subspace $\bm{H}^1(\Oref)$ with trace zero on
$\hat\Gamma^D\coloneqq \hat\Gamma_f^D\cup\hat\Gamma_s^D$ and
$\hat L\coloneqq L^2(\FLref)/\mathbb{R}$. Furthermore, let
$\vbref_D, \ubref_D\in\bm{H}^1(\Oref)$ be prolongations of the Dirichlet data
for the velocity and deformation, a right-hand side fluid force
$\fbref\in L^2(\FLref)$ be given. We define the \emph{stationary 
thermo-fluid-structure interaction problem} as follows.
Find $\vbref\in \{\Vbref + \vbref_D\}$, $\ubref\in
\{\Vbref+\ubref_D\}$, $\hat p\in \hat L$, and $\hat\theta\in\Vref$, such that
\begin{subequations}\label{eqn.ale-fsi}
  \begin{align}
    && \begin{multlined}[b]
    \big(\hat J\sigbref_{Bou}\hat F^{-T},\hat \nabla\phiref\big)_{\FLref}
    + \rho \big((\hat J\hat F^{-1} \vbref\cdot \hat\nabla) \vbref,\phiref\big)_{\FLref}
    - \big(\hat J \alpha_{\theta} (\hat\theta - \hat\theta_0) \fbref_f,\phiref\big)_{\FLref}\\
      + \big(\hat J\sigbref_R\hat F^{-T},\hat \nabla\phiref\big)_{\SOref}
    \end{multlined}
      &= \big(\rho_f \hat J\fbref,\phiref\big)_{\FLref} &\forall\phiref&\in\Vbref,\label{eqn.ale-fsi:momentum}\\
    && - \big(\vbref,\psiref\big)_{\SOref} +
      \big(\alpha_u \hat \nabla \ubref,\hat \nabla\psiref\big)_{\FLref}
      &=0 &\forall\psiref&\in\Vbref,\label{eqn.ale-fsi:ale}\\
    && \big(\widehat{\diver}\,(\hat J\hat F^{-1}\vbref_f),\xiref\big)_{\FLref}
      &=0&\forall\xiref&\in \hat L,\label{eqn.ale-fsi:mass}\\
    && \big(\hat J\hat\sigma_{\theta}\hat F^{-T},\hat\nabla\zetaref\big)_{\Oref}
      + \big(\hat J\hat F^{-1} \vbref\cdot \hat\nabla \hat\theta,\zetaref\big)_{\Oref}
      &= \big(\hat J f_{\theta},\zetaref\big)_{\Oref}&\forall\zetaref&\in\Vref,\label{eqn.ale-fsi:temp}
  \end{align}
\end{subequations}
with the harmonic mesh extension parameter $\alpha_u>0$, the stress
tensor in the solid as defined in \eqref{eqn.solid.cauchy-stress} and
$\hat J\sigbref_R\hat F^{-T} \coloneqq  \sigb_R$.
The thermal stress is $\sigma_\theta = k(\theta)\nabla\theta$ with
$\hat J\hat\sigma_\theta\hat F^{-T} \coloneqq  \sigma_\theta$.
The temperature dependent ALE fluid stress tensor $\sigbref_{Bou}$ is given by
\begin{equation*}
  \sigbref_{Bou} \coloneqq -\hat p_fI +\rho_f \nu(\hat\nabla \vbref_f \hat F^{-1}
  + \hat F^{-T}\hat\nabla \vbref_f^T),
\end{equation*}
with the kinematic viscosity $\nu>0$ and the fluid's density $\rho_f>0$.
\end{definition}

Let us comment on the above system in more detail. In \eqref{eqn.ale-fsi:momentum}, 
we combined the momentum equations of the fluid and the solid into one single equation.
This possible with variational-monolithic coupling
in which the interface conditions $\vbref_f = \vbref_s$ (Dirichlet) and
$\sigbref_{Bou}\nb = \sigbref_{R}\nb$ (Neumann) are fulfilled in an exact fashion
on the variational level. Moreover, the geometric condition
$\ubref_f = \ubref_s$ (Dirichlet) is fulfilled as well. The Dirichlet type conditions
are built into the function spaces as usual. The Neumann type condition
cancels out on the interface; see e.g., \cite[Section 3.3.3.5]{Wi20}.
In the second equation in \eqref{eqn.ale-fsi:ale}, the ALE mapping is realized
and for implementational reasons by using globally defined functions, we also work
with $\vbref = 0$ in $\Omega_s$. The third equation \eqref{eqn.ale-fsi:mass}
is the mass conservation of the fluid. The last equation \eqref{eqn.ale-fsi:temp}
is the weak form of the temperature equation.

\section{Numerical Tests}
\label{sec.numerical_tests}

In this section, we provide several numerical examples to validate and demonstrate the capabilities of the proposed algorithm. The numerical realisation is performed using Netgen/NGSolve\footnote{See also \url{https://ngsolve.org}}~\cite{Sch97,Sch14} and the add-on package ngsxfem~\cite{LHPvW21}.

\subsection{Numerical Approximation}
\label{sec.numerical_approx}

In total, our algorithm requires four different steps.

\begin{description}
    \item[Phase-field Fracture] A PFF problem must be computed both during initialization and in Step 4 of our algorithm. Specifically, we need to numerically solve the weak formulation of the phase-field fracture problem as described in \Cref{form_2_interface_b}. For this purpose, we use a finite element discretization, where both the phase-field and displacement field spaces are discretized using continuous piecewise linear finite elements. The pressure and temperature are provided as external parameters for this problem.

    To initialize each phase-field computation, we set $\varphi=0$ in $\Omega_f$ and $\varphi = 1$ else, and then solve \eqref{eqn.form_2_interface_b.b} without the coupling terms, applying a homogeneous Neumann boundary condition. This ensures that the initial condition satisfies the phase-field equation, preventing artificial strength at the crack tips.
    
    \item[COD computation] 
    In Step 1 of our algorithm, to compute the crack opening displacements (CODs) from the approximated phase-field function and displacement field, we use the unfitted finite element technology provided by ngsxfem to evaluate \eqref{eqn:cod-definition}  over an arbitrary line defined by a level set function. Notably, this level set does not need to be aligned with the mesh.

    \item[Geometry reconstruction] To construct the geometry from the computed CODs, we use Netgen's OCC (OpenCascade) interface to create a piecewise linear approximation of the interface, which is then meshed. Consequently, this step can also be viewed as an automated CAD model generation process. This is Step 2 of our algorithm.
    
    \item[Thermo-fluid-structure-interaction] 
    In Step 3 of our algorithm, we numerically solve the weak formulation of the thermo-fluid-structure interaction problem in ALE coordinates as defined in \Cref{def:fsi:ale:stationary}. We use the given mesh and discretize the spaces with inf-sup stable elements. Specifically, the velocity space is discretized using continuous piecewise quadratic elements, the pressure with continuous piecewise linear elements, and both the displacement and temperature with continuous piecewise quadratic elements.
\end{description}

\subsection{Example 1: \texorpdfstring{Convergence and}{} Temperature Sensitivity Study}
\label{sec.numerical_approx:subsec.ex1}

As a first example, we consider a basic test inspired by Sneddon's test~\cite{Sne46,SneddLow69}. Here, we do not consider the full model, but only the phase-field problem given in \Cref{form_2_interface_b}, to validate our model as used in Step 0 of our algorithm.

\paragraph{Set-up}
The set-up for this is as follows. We consider $\Omega=(-2,2)^2$ and an the initial crack (where $\varphi^0=0$) is given in  $(-0.2, 0.2)\times(-\hcr,\hcr)$. The material parameters are $E=1$, $\nu_s=0.3$, $G_c=1$. The Lamé parameters are then obtained by $\mu = \frac{E}{2(1+\nu_s)}$ and $\lambda = \frac{\nu_s E}{(1+\nu_s)(1-2\nu_s)}$. 
The pressure is $p=4\times10^{-2}$ and the reference pressure is $p_0=0$. Both the temperature and reference temperature are $\theta=\theta_0=0$. The discretizations parameters are $\kappa=10^{-10}$, $\varepsilon= 2 h$, where $h$ is the local mesh size, and $\gamma = 100 h^{-2}$. 

\begin{figure}
  \centering
  \includegraphics{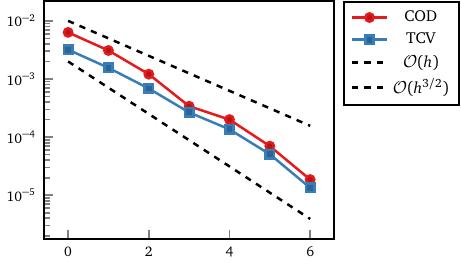}
  \caption{Numerical convergence of the crack opening displacement in the center of the crack and total crack volume in \hyperref[sec.numerical_approx:subsec.ex1]{Example 1}.}
  \label{fig.ex1.convergence}
\end{figure}

\paragraph{Convergence Results}
We consider this set up over a series of eight meshes, constructed such that the local mesh size of the crack $\hcr = \hmax / 100$. On the coarsest mesh level, we take $\hmax=1.28$. Each subsequent mesh is constructed by halving both mesh parameters. We compute the crack opening displacement (COD) in the center of the crack and the total crack volume (TCV) and compare the results with those obtained the finest mesh as the reference solution. This reference solution was computed with a total of $2.3\times 10^6$ degrees of freedom, and therefore close to the limit of the direct solver used to solve the resulting linear systems within each Newton-Step. The results for this can be seen in \Cref{fig.ex1.convergence}. We see that both the COD and TCV converge with a rate between 1 and 1.5.

\paragraph{Temperature effects}
To study the effect of the temperature difference on the crack opening displacement, we consider the previous set up on mesh level six with temperatures $\theta\in\{240, 160, 80, 0, -80, -160, -240\}$. Note that with our chosen material parameters we obtain $K_{dr} = \frac56$ and $\alpha_{\theta} - 3\alpha_{\theta} K_{dr} = - 1.5\times 10^{-5}$. Therefore, as we consider both constant pressure and temperature in this example, an increase in temperature of $\frac23 \times 10^5$ corresponds to a pressure decrease of $1$.

\begin{figure}
  \centering
  \includegraphics{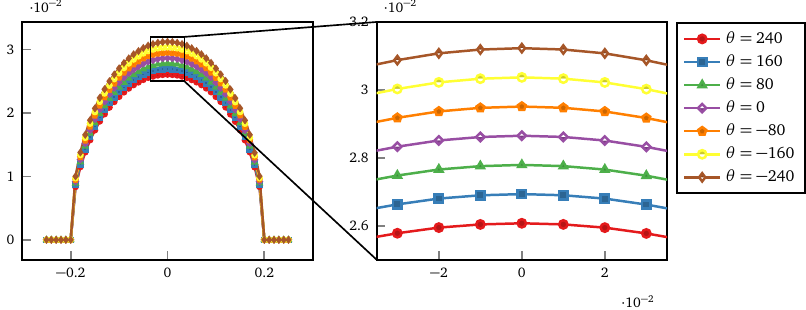}
  \caption{Crack Opening displacement dependent on temperature in \hyperref[sec.numerical_approx:subsec.ex1]{Example 1}.}
  \label{fig:ex1.change-temp}
\end{figure}

We present the results for the change of temperature in \Cref{fig:ex1.change-temp}. As expected, an increase in temperature relative to the reference temperature causes the crack width to shrink, while a decrease in temperature causes the crack to open further~\cite{morrow2001permeability,yasuhara2006evolution,hardin1982measuring,rutqvist2008analysis}.

\subsection{Example 2: Fully Coupled Stationary Test Case}
\label{sec.numerical_approx:subsec.ex2}

As our second example, we consider the fully coupled algorithm in a stationary setting.

\paragraph{Set-up}
The basic set up is the same as in \hyperref[sec.numerical_approx:subsec.ex1]{Example 1}. 
That is, we consider $\Omega=(-2,2)^2$, the initial crack is $(-0.2, 0.2)\times(-\hcr,\hcr)$, the material parameters are $E=1$, $\nu_s=0.3$, $G_c=1$, the pressure is $p=p^0=4\times10^{-2}$, the reference pressure is $p_0=0$, and both initial temperature and reference temperature are $\theta=\theta_0=0$. For the thermo-fluid-structure interaction problem we have $\mathbf{\kappa}_f = 0.01$, $\mathbf{\kappa}_s=1.0$ and the force per unit mass is $1\eb_z$. The fluid is driven by the body force $\fbref = 0.2 \exp(-1000\vert \xb \vert^2)\eb_x$ and the temperature is driven by the external forcing term $\hat{f}_\theta = 100 \exp(-10\vert \xb \vert^2)$. The maximal mesh size chosen is $h=0.16$ and the crack has a local mesh size of $\hcr = 0.0016$. This corresponds to mesh level 3 in the previous example. The remaining discretisations parameters are as in \hyperref[sec.numerical_approx:subsec.ex1]{Example 1}.

In Step 4, the temperature from the TFSI problem is then used as the temperature in the PFF problem, i.e., $\theta^n = \theta^{n, TFSI}$. As the pressure in the TFSI problem is normalized to be mean zero, we add this to the initial temperature as the driving pressure in the PFF problem, i.e., $p^n= p^0 + p^{n, TFSI}$.

\begin{figure}
  \centering
  \includegraphics{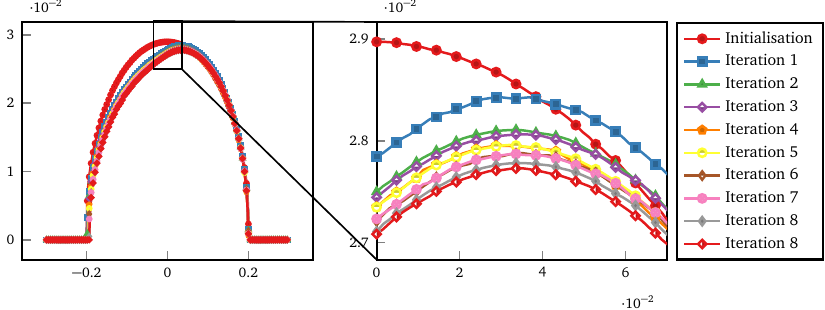}
  \caption{Crack Opening displacement after each TFSI and PFF iteration in \hyperref[sec.numerical_approx:subsec.ex2]{Example 2}. Marks at every third computed COD.}
  \label{fig:exmpl2.cod}
\end{figure}

\begin{figure}
  \centering
  \begin{subfigure}[t]{0.325\textwidth}
    \centering
    \includegraphics[width=.95\linewidth]{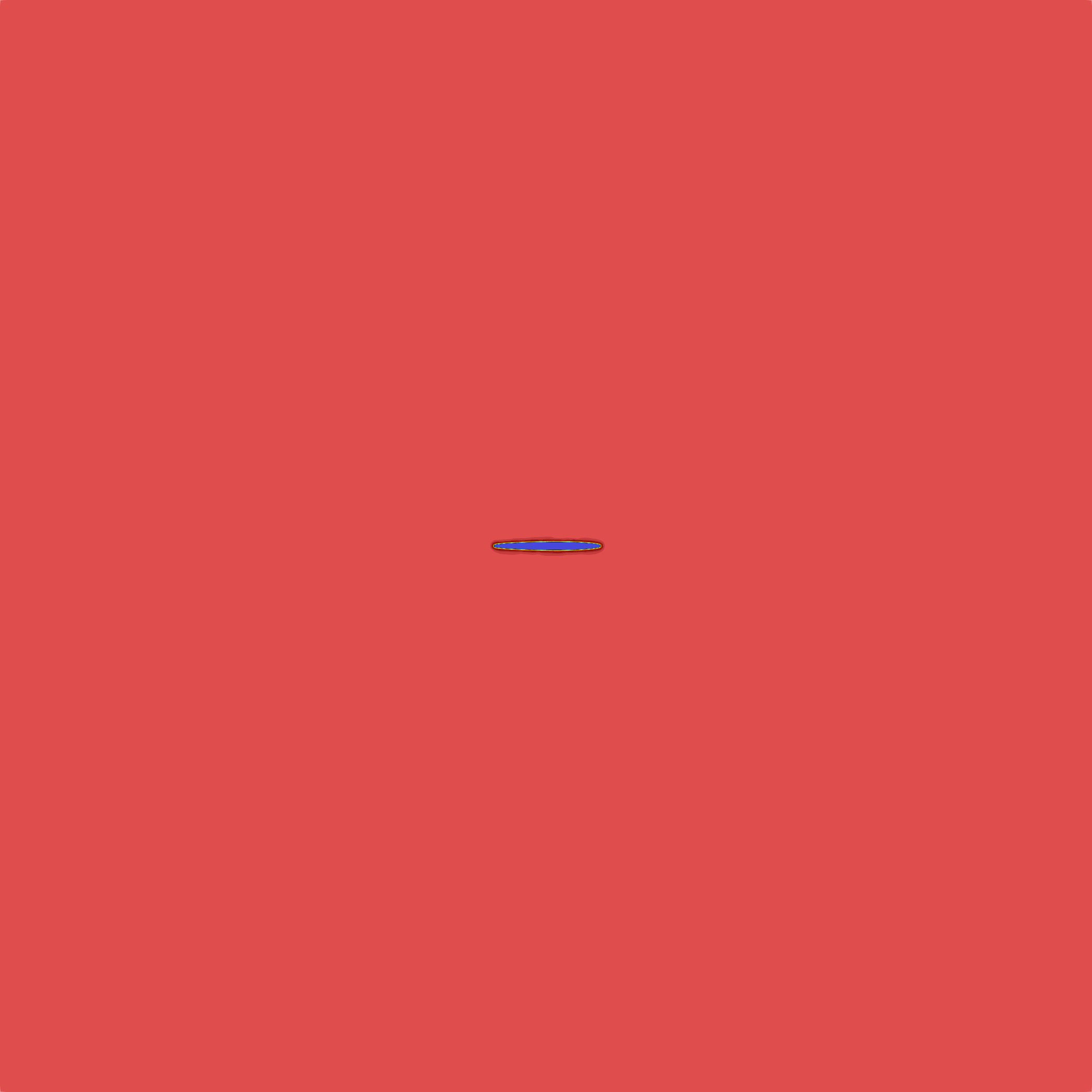}
    \vspace*{18.5pt}
    \caption{Phase field}
  \end{subfigure}
  \hfill
  \begin{subfigure}[t]{0.325\textwidth}
    \centering
    \includegraphics[width=.95\linewidth]{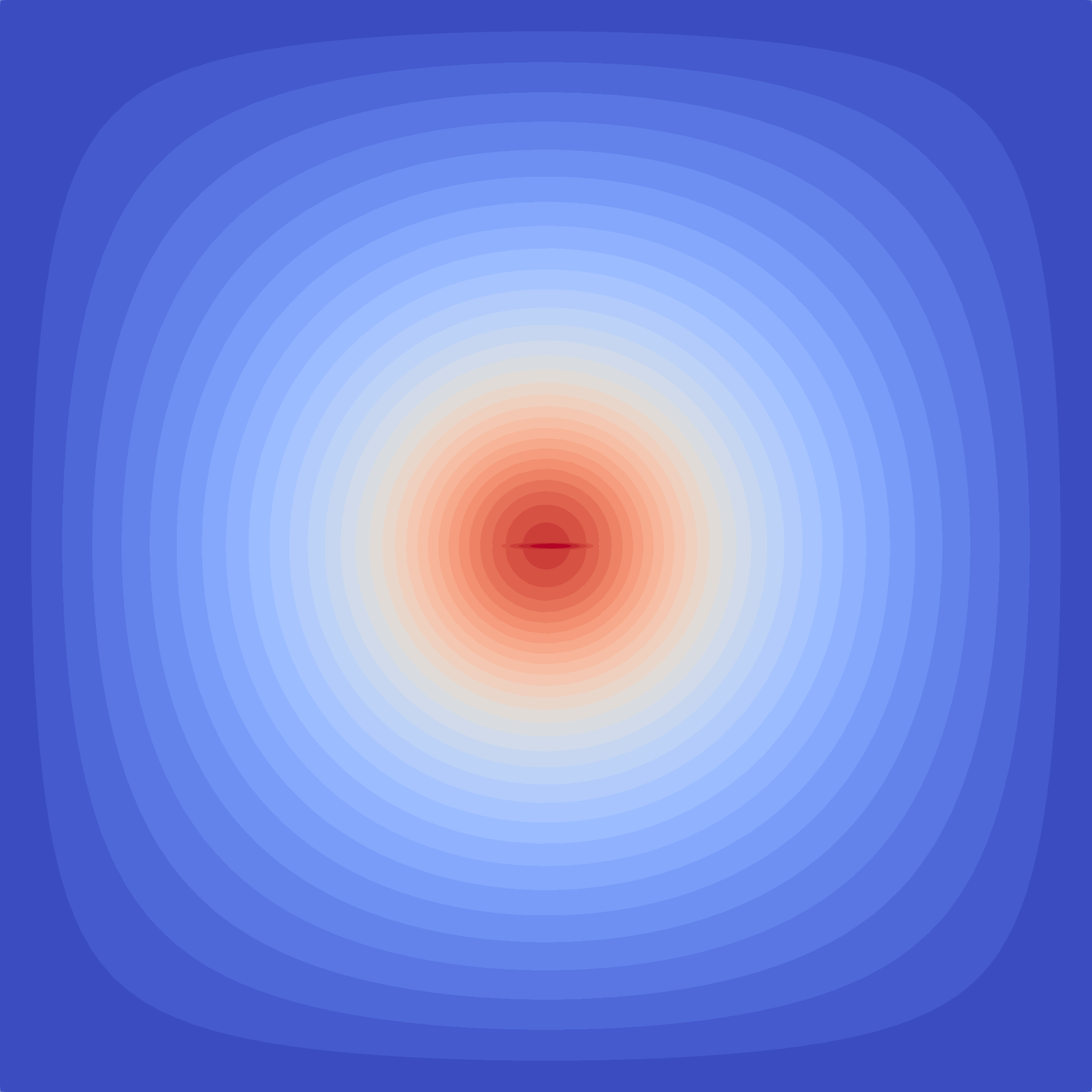}\\[5pt]
    \includegraphics[height=12.5pt]{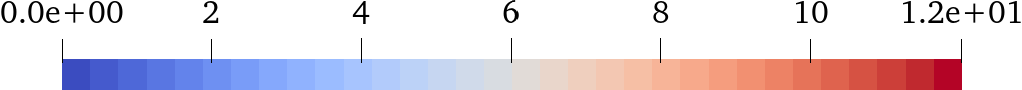}
    \caption{Temperature}
  \end{subfigure}
  \hfill
  \begin{subfigure}[t]{0.325\textwidth}
    \centering
    \includegraphics[width=.95\linewidth]{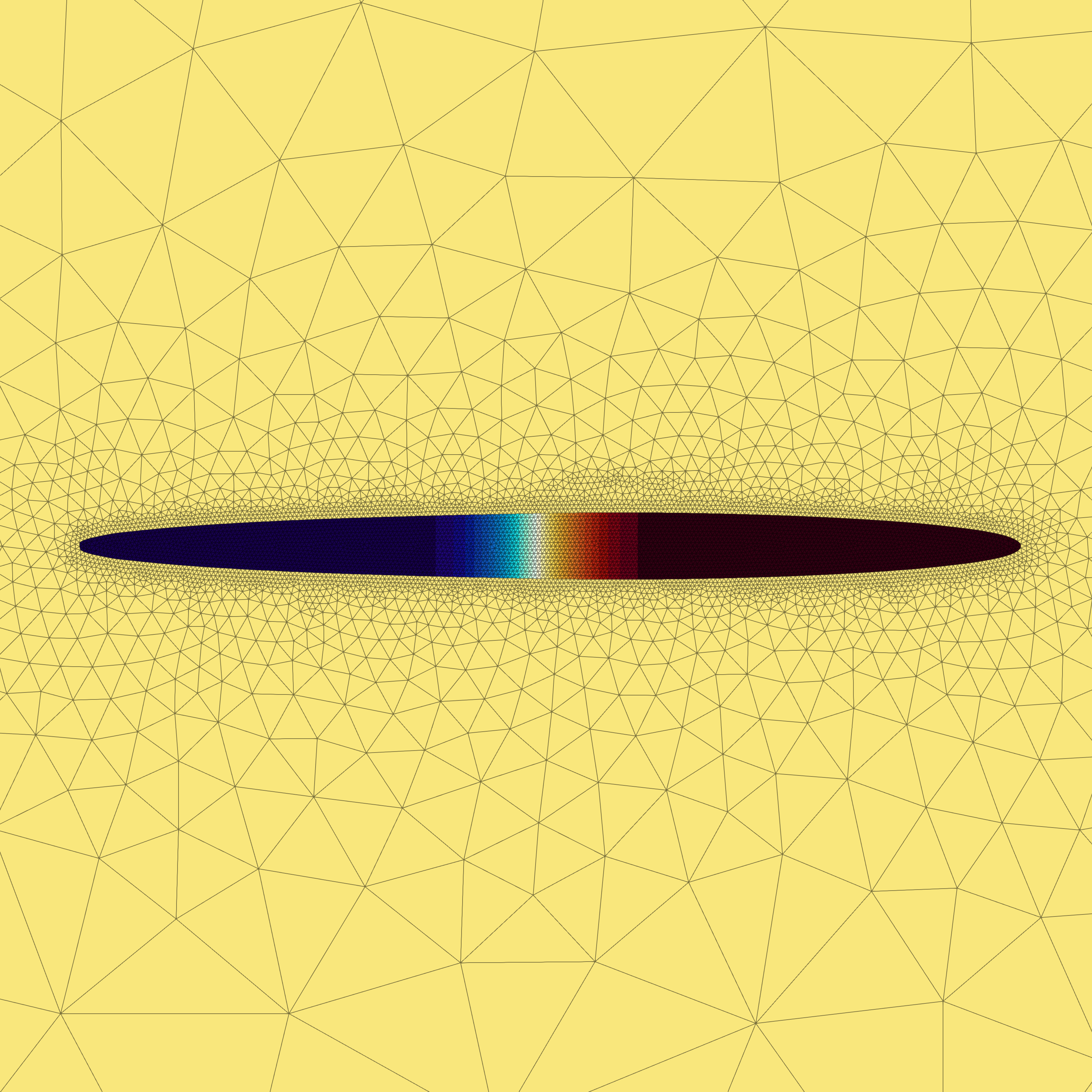}\\[5pt]
    \includegraphics[height=12.5pt]{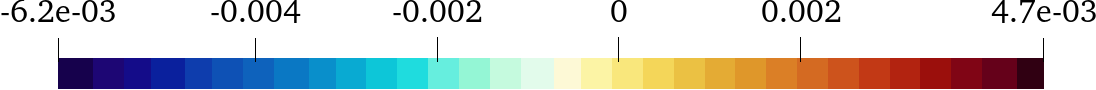}
    \caption{Pressure and mesh}
  \end{subfigure}
  \caption{Phase-field, temperature and pressure (with mesh) after the last iteration in \hyperref[sec.numerical_approx:subsec.ex2]{Example 2}.}
  \label{fig.ex2:images}
\end{figure}

\paragraph{Results}
The resulting crack opening displacements for eight iterations between the phase-field-fracture and thermo-fluid-structure-interaction problem can be seen in \Cref{fig:exmpl2.cod}, and the resulting total crack volume for every iteration in \Cref{tab:ex2.tcv}. A visualization of the phase-field, temperature and pressure after the final iteration can be seen in \Cref{fig.ex2:images}. As the FSI temperature is positive (relative to the reference temperature), we have the expected reduction in the fluid/crack volume as seen through the COD plot and TCV values. Furthermore, the FSI pressure is negative in the left half of the domain and positive in the right half of the crack, resulting in the observed skewness of the crack. Furthermore, we see that after eight iterations, there is very little change in the COD.

\begin{table}
  \centering
  \caption{Total crack volume after each TFSI and PFF iteration in \hyperref[sec.numerical_approx:subsec.ex2]{Example 2}.}
  \label{tab:ex2.tcv}
  \begin{tabular}{rc@{\hspace*{5pt}}c@{\hspace*{5pt}}c@{\hspace*{5pt}}c@{\hspace*{5pt}}c@{\hspace*{5pt}}c@{\hspace*{5pt}}c@{\hspace*{5pt}}c@{\hspace*{5pt}}c@{\hspace*{5pt}}c}
    \toprule
    Iteration & Init. & 1&2&3&4&5&6&7&8&9\\  
    \midrule
    TCV  & 0.00926 & 0.00883 & 0.00869 & 0.00868 & 0.00859 & 0.00859 & 0.00854 & 0.00854 & 0.00847 & 0.00845\\
    \bottomrule
  \end{tabular}
\end{table}

\paragraph{Domain reconstruction}
For fine-meshes and and near the singularity of the crack tip, the computation of the COD becomes numerically unstable. Due to the increasing number of necessary CODs computed with decreasing mesh size, this causes a rough boundary of the approximated crack and with each iteration of the coupled loop, this effect increases. To smooth out the crack boundary for the FSI (and subsequent PFF) computations, we process the COD data before the domain is reconstructed. To this end, we use the \texttt{numpy} function $\texttt{numpy.polynomial.chebyshev.chebfit}$ to compute a least-squares polynomial approximation of the crack boundary. The resulting polynomial values at the $\pazocal{O}(h^{-1})$ points where the COD was originally computed, plus the roots of this polynomial, are then used to define the crack boundary.

\begin{figure}
  \centering
  \includegraphics{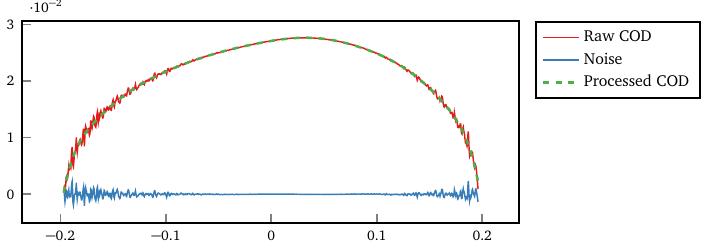}
  \caption{Smoothed and original COD data for \hyperref[sec.numerical_approx:subsec.ex2]{Example 2} with $\hmax=0.04$ after four iterations of the coupled loop computed on the rough COD reconstructed domain.}
  \label{fig.ex2.cod-smoothing}
\end{figure}

In \Cref{fig.ex2.cod-smoothing}, we illustrate the process on the COD data resulting after four iterations of our coupling loop where the smoothing was not applied using a mesh with $\hmax=0.04$ and $\hcr=0.0004$. Here, we see the visible oscillations towards the crack tip are smoothed out effectively while the shape of the crack is maintained.

\subsection{Example 3: Propagating Crack}
\label{sec.numerical_approx:subsec.ex3}

We test our novel phase-field model in the case of a propagating crack in two situations: A spatially constant and increasing pressure and finally with a time-dependent temperature and pressure resulting from a TFSI problem. The first example aims to study pressure-driven crack propagation for our novel interface phase-field approach, while the latter aims to study the effect of the temperature coupling to the crack's propagation.

\subsubsection{Example 3a: Pressure driven propagation without TFSI coupling}
\label{sec.numerical_approx:subsec.ex3a}

To test the present interface phase-field model regarding crack propagation, we consider a series of loading steps (pseudo time-steps) and apply an increasing pressure in each iteration, rather than taking the pressure and temperature from a TFSI problem. As the temperature and pressure enter the phase-field model through the same interface integrals, it is sufficient to just consider the pressure in this example.

\paragraph{Set-up}
The basic set-up is as before. We have, $\Omega=(-2,2)^2$, the initial crack is $(-0.2, 0.2)\times(-\hcr,\hcr)$, material parameters are $E=1$, $\nu_s=0.3$, $G_c=1$, the initial pressure is $p=4\times10^{-2}$ and the temperature is $\theta=0$. The pressure in each iteration of our loop is chosen as $p^n = 4\times 10^{-2} + n \times 10^{-4}$ and, we consider a total of $100$ loading steps.

We consider a series of three meshes. The coarsest mesh is constructed with $\hmax = 0.5$ and $\hcr = 0.005$. Each subsequent mesh is constructed by halving both mesh parameters. The penalization parameter is again $\gamma = 100 h^{-2}$, but the phase-field regularization parameter is fixed to $\varepsilon=0.01$ for all meshes. The latter corresponds to the previous choice on the coarsest mesh in this example.

\paragraph{Results}

\begin{figure}
  \centering
  \includegraphics{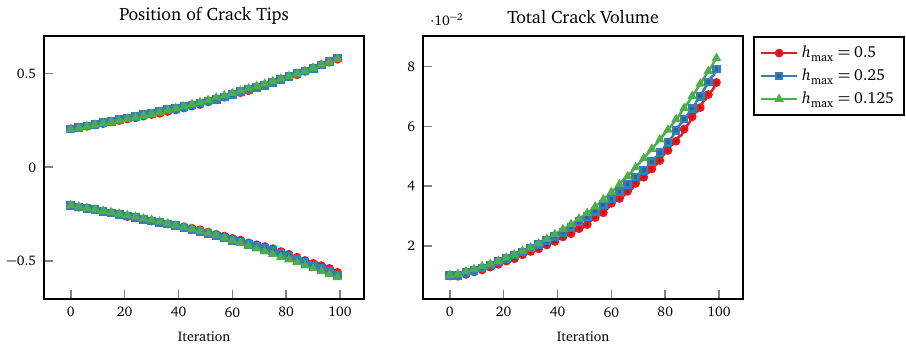}
  \caption{Horizontal position of crack tips and total crack volume resulting from a spatially constant, increasing pressure in \hyperref[sec.numerical_approx:subsec.ex3a]{Example 3a}. Marks at every third iteration.}
  \label{fig.ex3a:crack_tips_tcv}
\end{figure}

\begin{figure}
  \centering
  \includegraphics[width=0.19\linewidth]{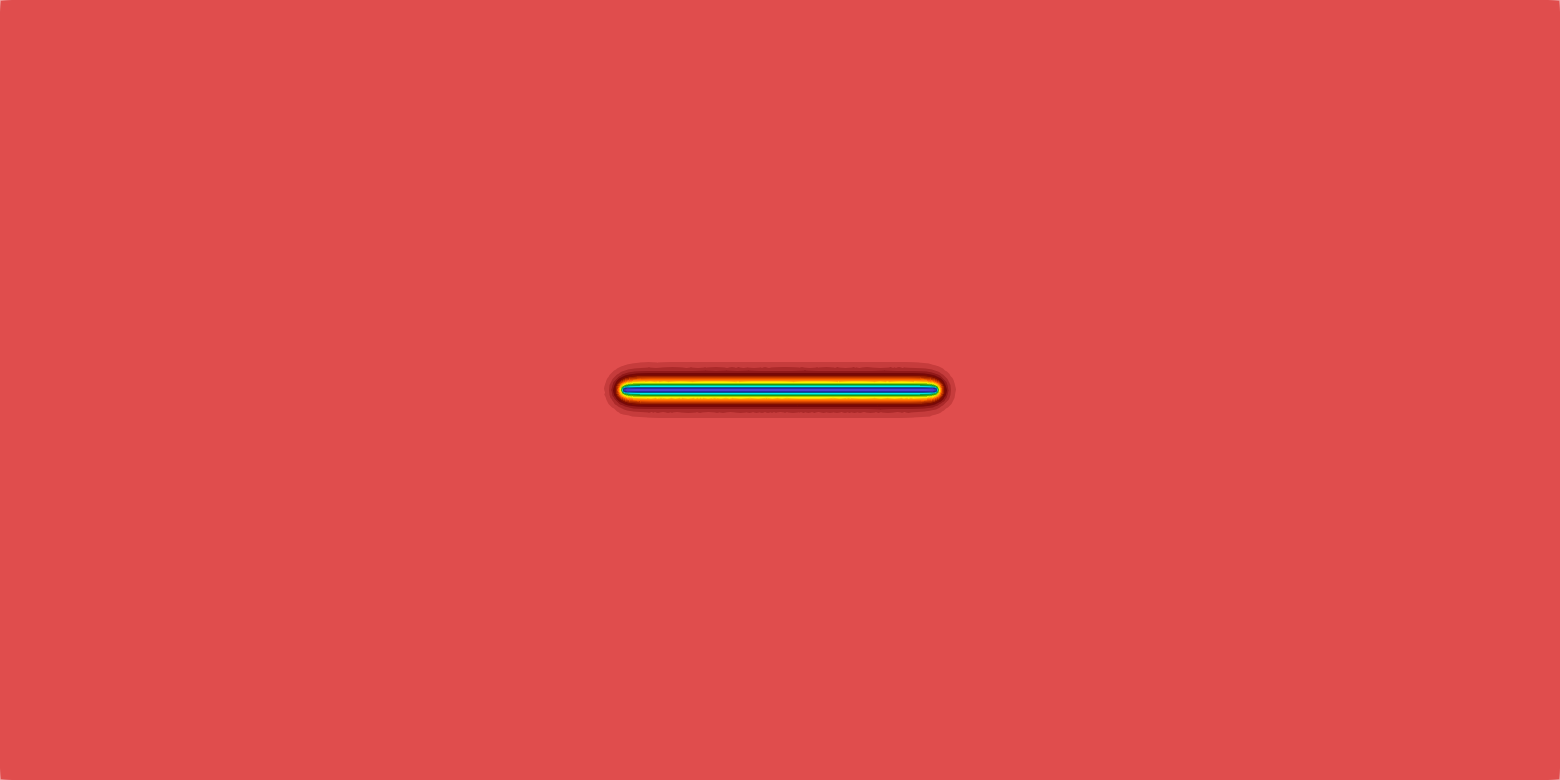}
  \includegraphics[width=0.19\linewidth]{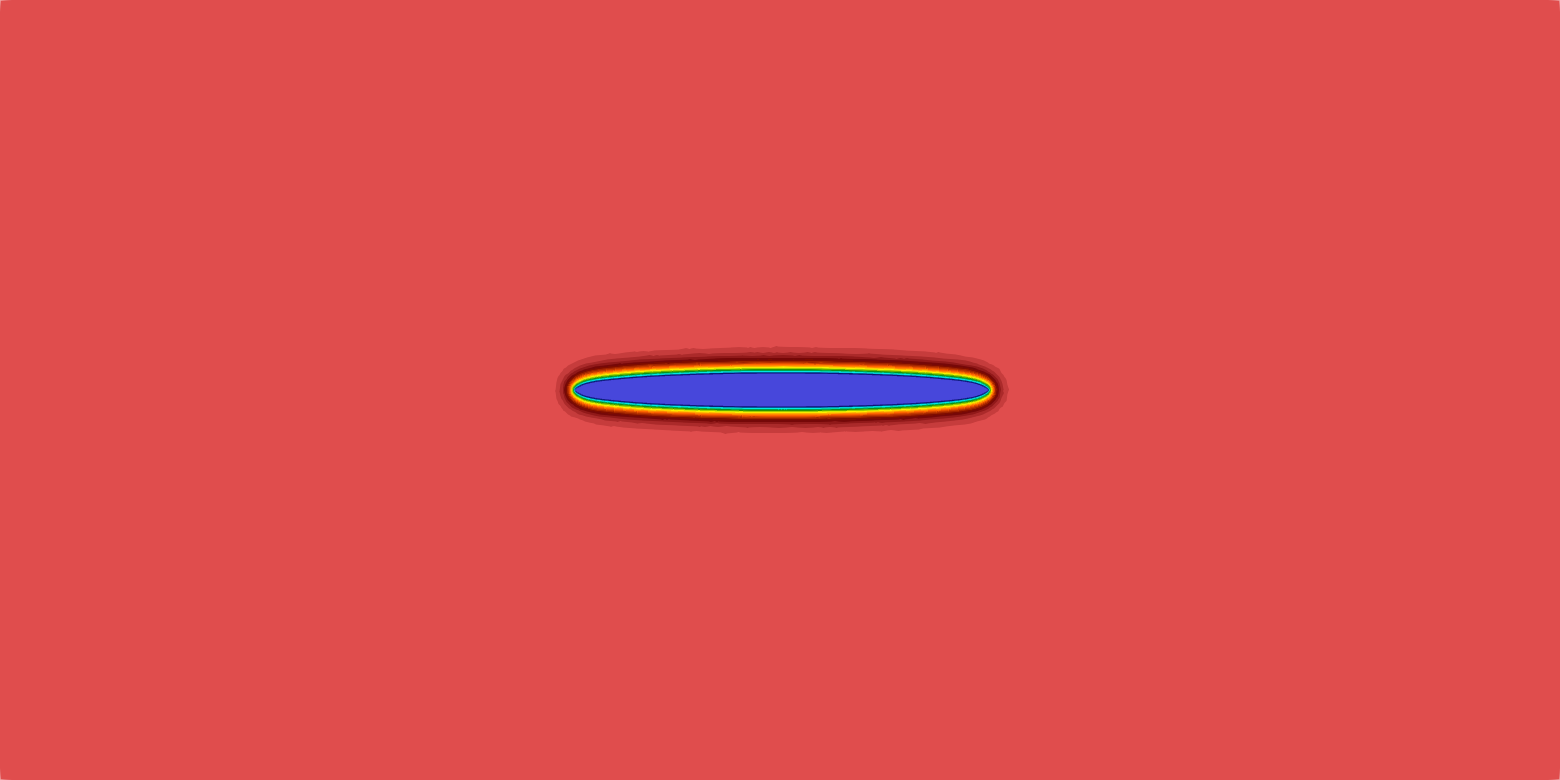}
  \includegraphics[width=0.19\linewidth]{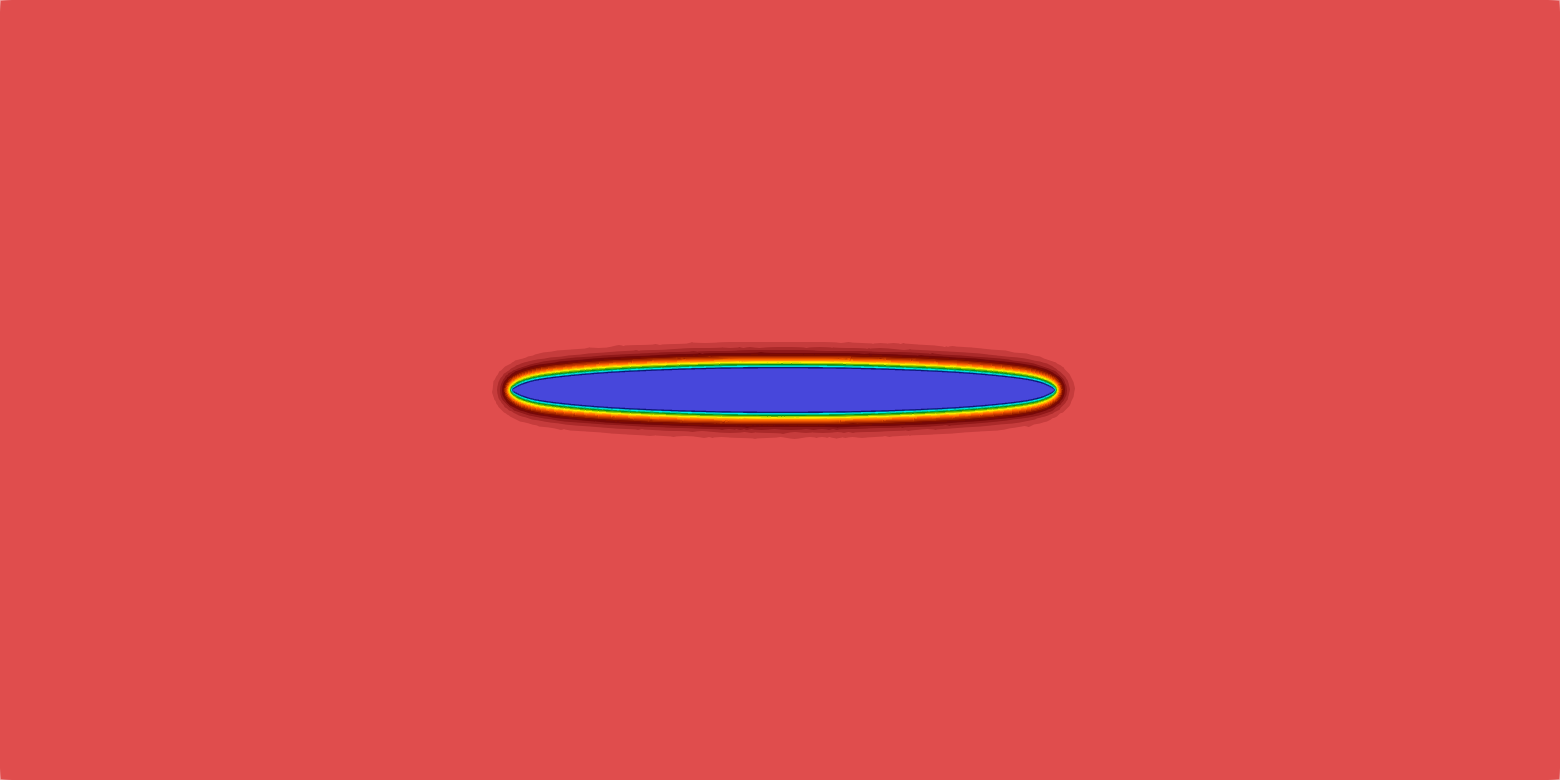}
  \includegraphics[width=0.19\linewidth]{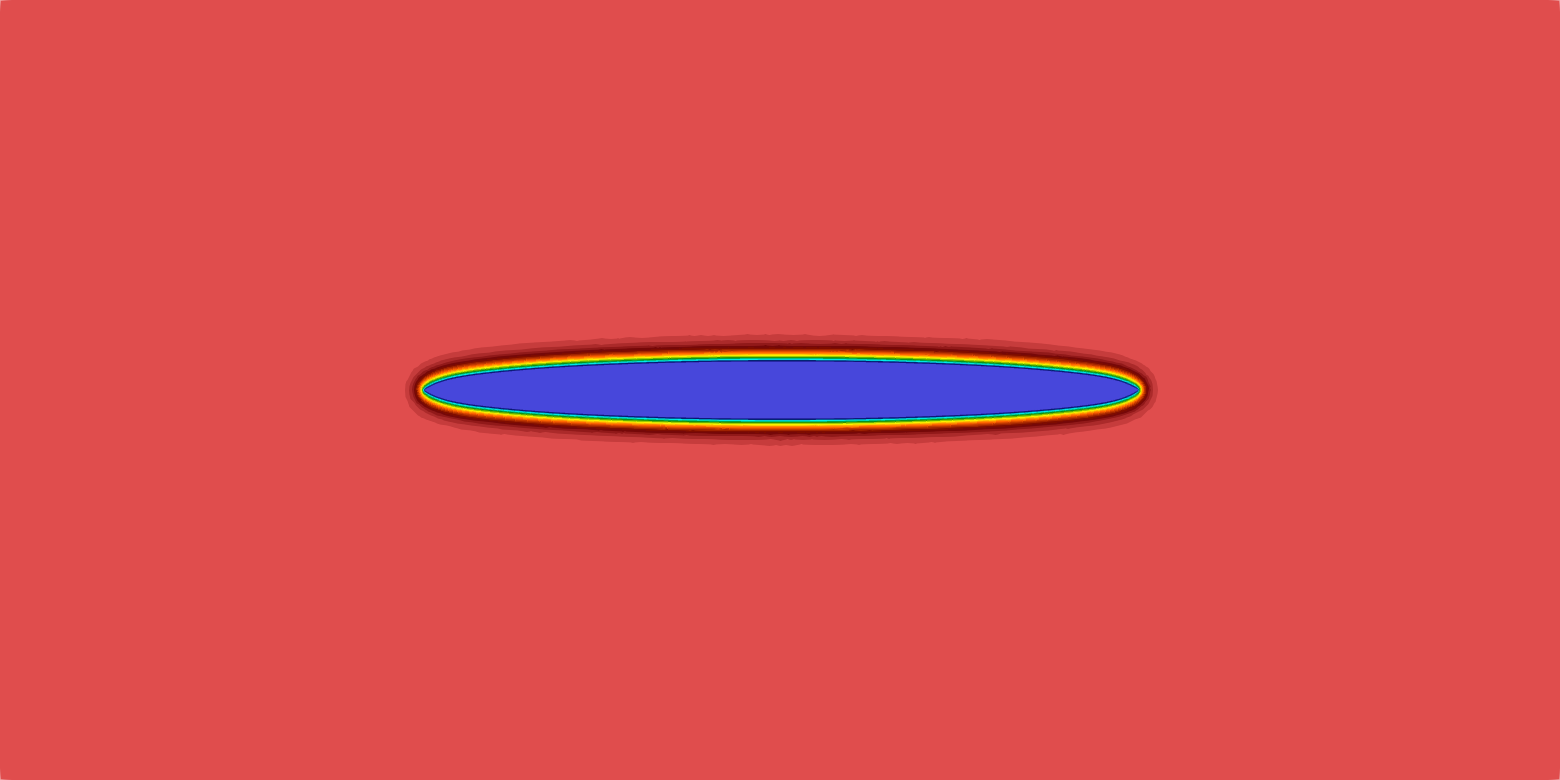}
  \includegraphics[width=0.19\linewidth]{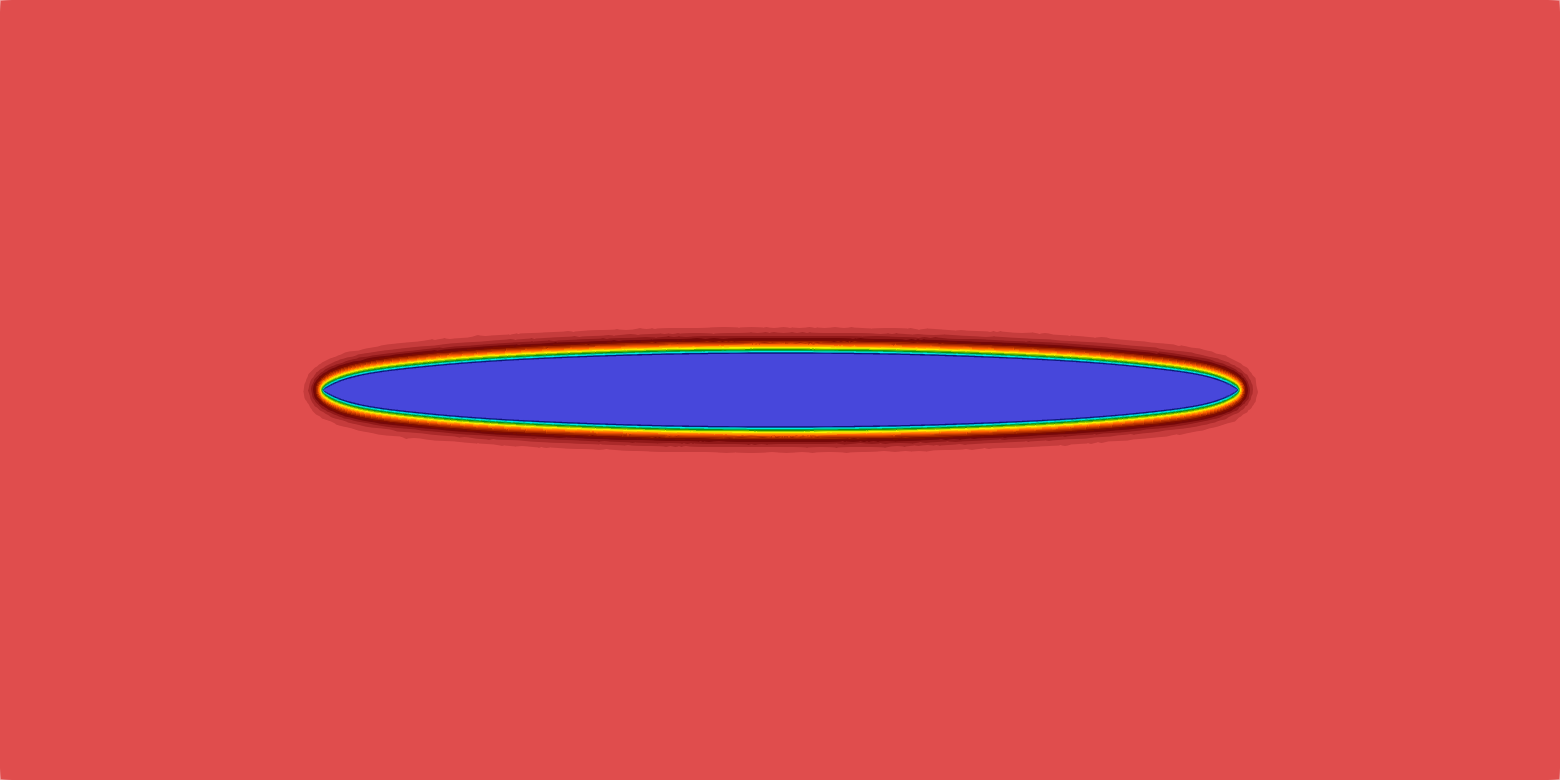}
  \caption{Phase field in the sub-domain $(-1, 1)\times(-0.5, 0.5)$ form the initialization and in iteration steps 25, 50, 75, 100 in \hyperref[sec.numerical_approx:subsec.ex3a]{Example 3a}. Computed with $\hmax=0.125$.}
  \label{fig.ex3a:phase-field}
\end{figure}

The resulting horizontal position of the crack tips and the total crack volume from each iteration are shown in \Cref{fig.ex3a:crack_tips_tcv}. We can see the consistent propagation of the crack over each of the four meshes, thereby showing that crack propagation is also feasible for our phase-field model. In \Cref{fig.ex3a:phase-field}, we show the phase-field at initialization and for four iterations. Here, we see that the phase-field is well behaved and that the crack increases both in the horizontal and vertical directions.

\subsubsection{Example 3b: Propagation with TFSI coupling}
\label{sec.numerical_approx:subsec.ex3b}

In this example, we consider the fully coupled scheme. 

\paragraph{Set-up}
The initialization step is as in \Cref{sec.numerical_approx:subsec.ex3a}. In the reconstructed domain, the TFSI problem is driven by the fluid and temperature forcing terms
\begin{equation*}
\fbref = 0.02 \exp( -400(x-0.1)^2 )\eb_x\quad\text{and}\quad \hat{f}_\theta= -800\exp(-\vert\xb\vert^2).
\end{equation*}
The heat conductivity parameters are $\mathbf{\kappa}_f = 0.005$, $\mathbf{\kappa}_s=1.0$ and the force per unit mass is $0.5\eb_z$. The remaining material parameters are as in \hyperref[sec.numerical_approx:subsec.ex3a]{Example 3a}.

In each iteration of our coupling loop, the phase-field is then driven by the pressure $p^n = 4\times 10^{-2} + 5n\times10^{-5} + p^{n,TFSI}$ and temperature $\theta^n=\theta^{n, TFSI}$. We again consider 100 loading steps. The remaining phase-field parameters are again chosen as in \hyperref[sec.numerical_approx:subsec.ex3a]{Example 3a}.

\begin{figure}
  \centering
  \begin{subfigure}[t]{.49\textwidth}
    \centering
    \includegraphics[width=.6\textwidth]{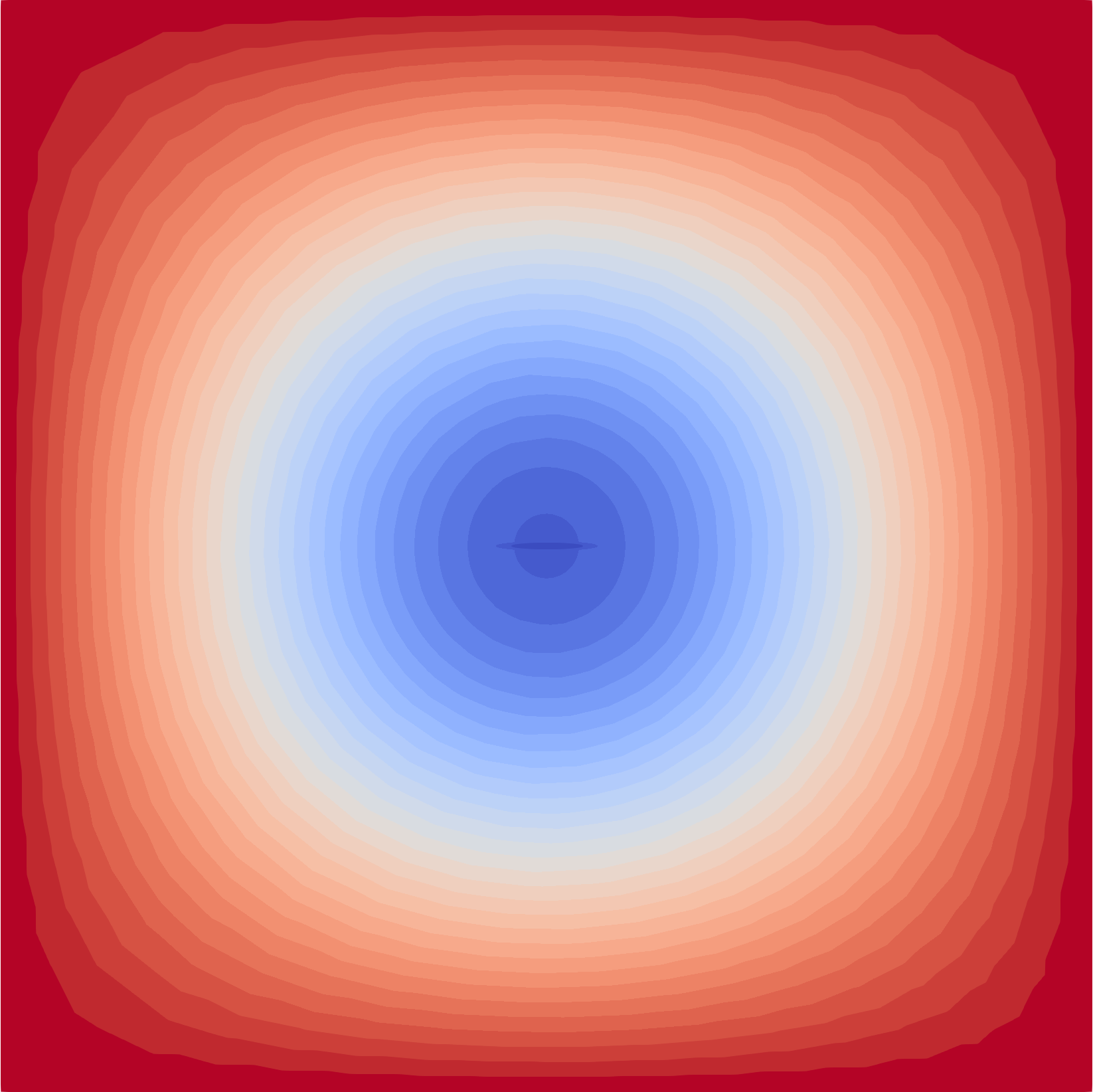}\\[5pt]
    \includegraphics[height=12.5pt]{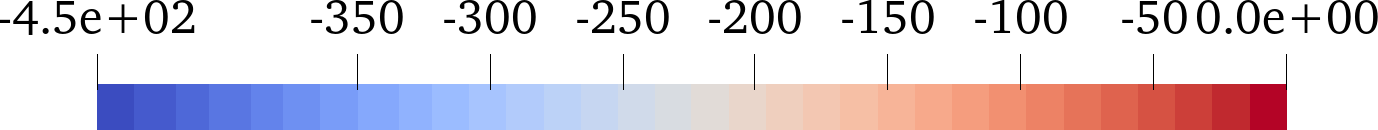}
    \caption{Temperature}
  \end{subfigure}
  \begin{subfigure}[t]{.49\textwidth}
    \centering
    \includegraphics[width=.6\textwidth]{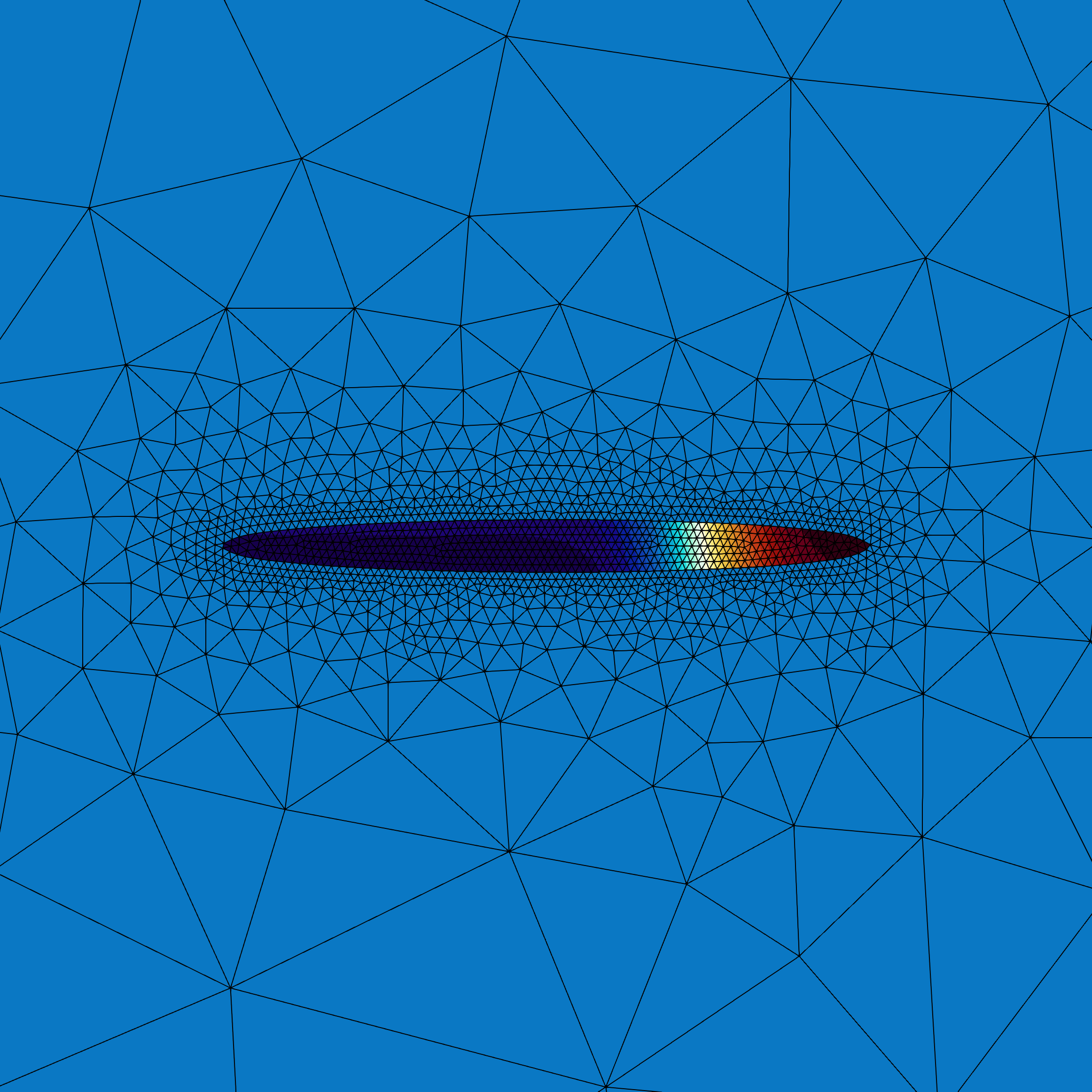}\\[5pt]
    \includegraphics[height=12.5pt]{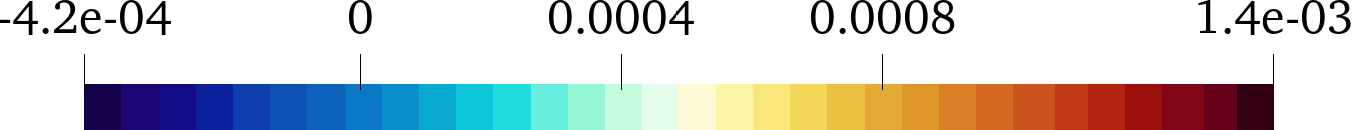}
    \caption{Pressure and mesh}
  \end{subfigure}
  \caption{TFSI temperature and pressure in the first iteration in \hyperref[sec.numerical_approx:subsec.ex3b]{Example 3b}. Computed on the mesh with $\hmax=0.5$.}
  \label{fig.ex3b:temp+pre}
\end{figure}

\paragraph{Results}
We consider a series of three meshes. First, we see in \Cref{fig.ex3b:temp+pre} the resulting temperature and pressure from the TFSI computation in the first iteration. Here we see that the temperature is negative and that the TFSI pressure is positive at the right tip of the crack and negative at the left tip of the crack. Due to the cooling of the medium and the higher pressure at the right tip of the crack, we expect the fracture to grow faster than in \hyperref[sec.numerical_approx:subsec.ex3a]{Example 3a}. Furthermore, we expect it to grow faster towards the left than the right.

\begin{figure}
  \centering
  \includegraphics{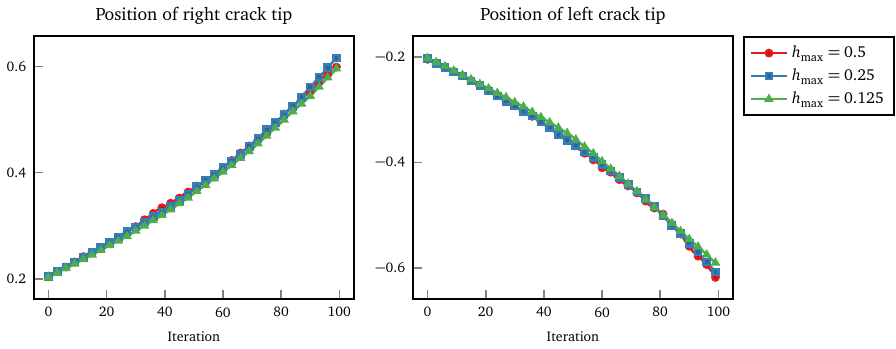}
  \caption{Position of the left and right crack tips in \hyperref[sec.numerical_approx:subsec.ex3b]{Example 3b}. Mark at every third iteration.}
  \label{fig.ex3b:crack-tips}
\end{figure}

In \Cref{fig.ex3b:crack-tips}, we see the position of the left and right tips of the crack for each considered mesh. First we note, that the crack does indeed grow faster than in \hyperref[sec.numerical_approx:subsec.ex3a]{Example 3a} and that the results are consistent over the series of meshes. Furthermore, we see that on finer meshes, the trajectory becomes smoother. After twenty iterations, the left and right positions of the crack tips are $-0.2591, 0.2656$, $-0.2609, 0.2651$ and $-0.2571, 0.2598$ for the three meshes, respectively. We, therefore, initially observe the expected faster growth toward the right than the left. However, this difference is not very large and after 100 iterations, the left tip has moved further than the right in some cases. We attribute this to numerical error, since this inconsistency occurs earlier on coarser meshes, which we observe to be less stable in \Cref{fig.ex3b:crack-tips}.

\subsection{Example 4: Two orthogonal cracks with TFSI coupling}
\label{sec.numerical_approx:subsec.ex4}

In this final example, we consider two orthogonal, connected cracks, illustrating that this approach is applicable to multiple cracks. While this situation is not challenging for phase-field computations, it is more involved with respect to the geometry reconstruction.

To study the effects of including the temperature in our model, we compare two cases here. First, we only couple the TFSI pressure back to the phase-field computation, and denote the resulting phase-field fracture deformation by $\ub_h^{p}$. Secondly, we couple both the temperature and pressure back to the phase-field computation, and denote the resulting phase-field fracture deformation by $\ub_h^{p, \theta}$.

The following example is an extension of the example presented in \cite[Section~5.4]{WaWi23_CMAME}.

\paragraph{Set-up}
The background domain is $\O=(-2,2)^2$, and as before, we consider homogeneous Dirichlet boundary conditions for the displacement and homogeneous Neumann conditions for the phase-field. The initial phase field is given by a flipped "T", i.e.,
\begin{equation*}
  (-0.2, 0.2)\times(-h, h)\cup(0.2-h, 0.2+h)\times(-0.2, 0.2).
\end{equation*}
As before, the material parameters are as above $E=1$, $\nu_s=0.3$, and $G_c=1$, while the initialization pressure and temperature are $p=4\times 10^{-2}$ and $\theta=0$, respectively.

The thermo-fluid-structure interaction problem in the reconstructed domain is then driven by the forcing terms
\begin{equation*}
\fbref = 0.2\exp(-1000( (x - 0.2)^2 + y^2))\eb_x \quad\text{and}\quad \hat{f}_\theta= 100 \exp(-10 (x-0.2)^2 -5 y^2).
\end{equation*}
The reference temperature is set as $\theta_0 = 0$, the heat conductivity parameters are $\mathbf{\kappa}_f = 0.01$, $\mathbf{\kappa}_s=1.0$ and the force per unit mass is $1\eb_y$.

\paragraph{Results}

\begin{figure}
  \centering
  \includegraphics[width=.3\textwidth]{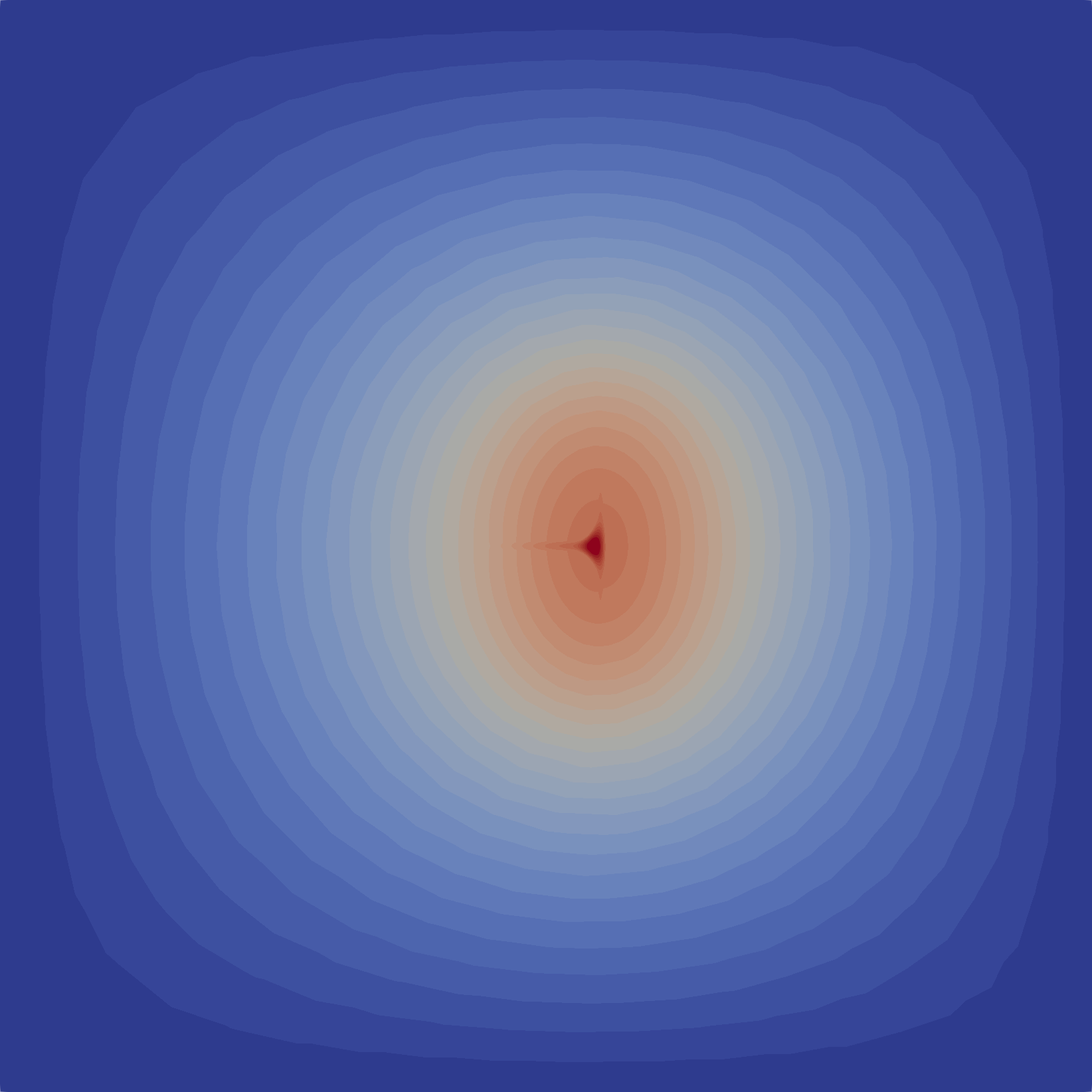}
  \includegraphics[width=.3\textwidth]{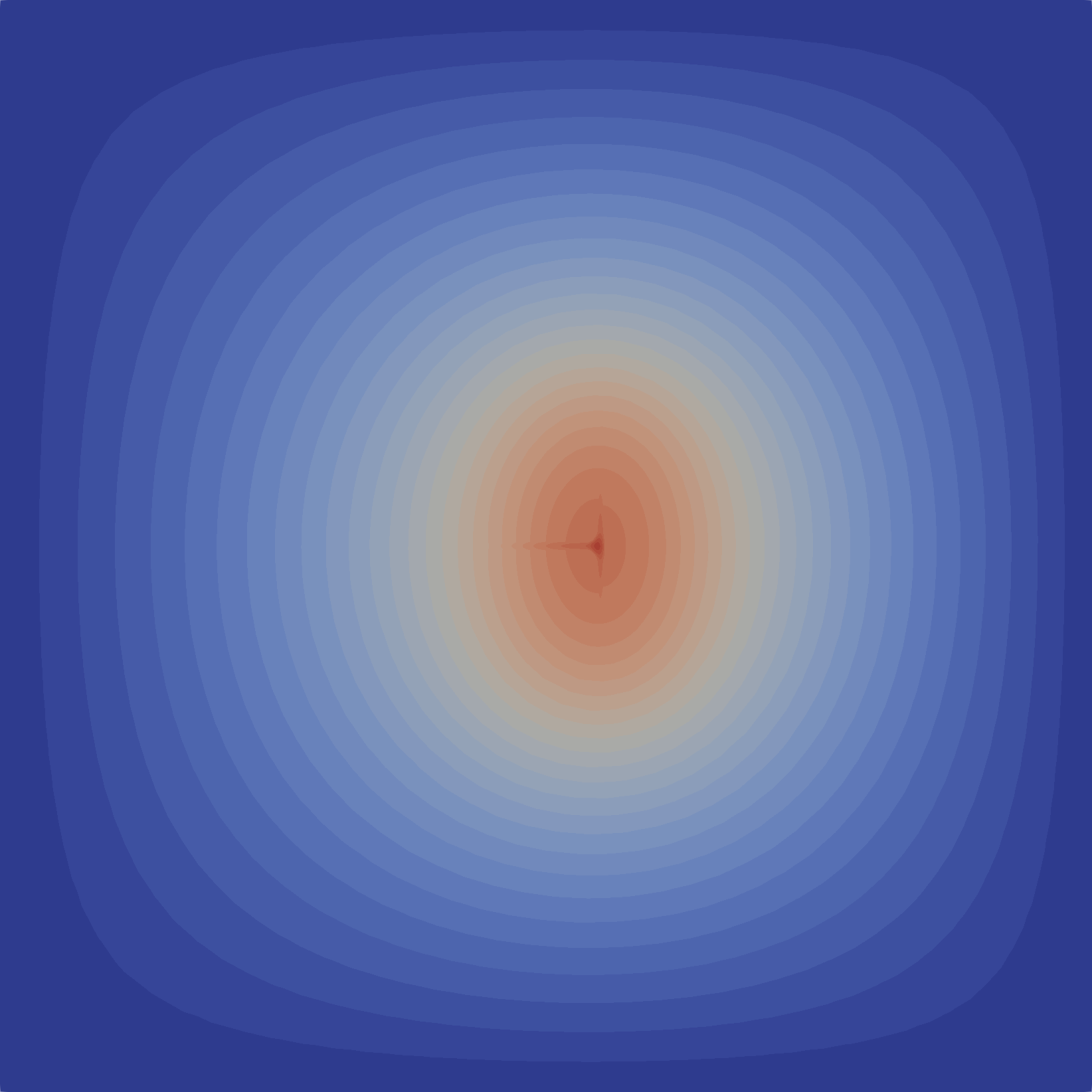}
  \includegraphics[width=.3\textwidth]{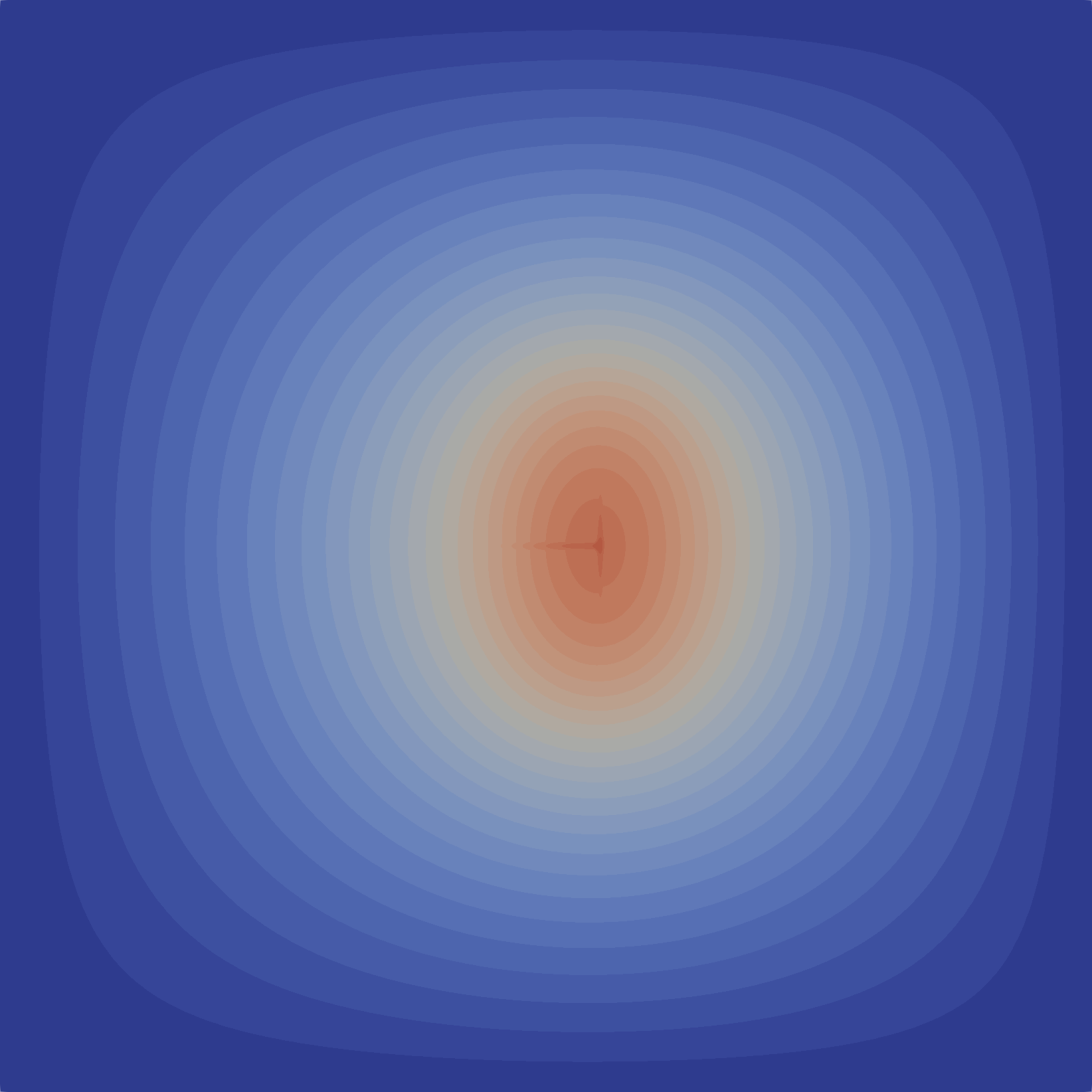}\\[5pt]
  \includegraphics[height=12.5pt]{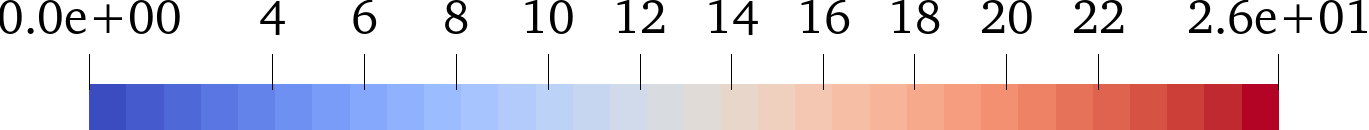}
  \caption{TFSI Temperature computed in the reconstructed domain in \hyperref[sec.numerical_approx:subsec.ex4]{Example 4}. Computed on meshes with $\hmax=0.5, 0.25$ and $0.125$, respectively.}
  \label{fig:ex4.tfsi-temp}
\end{figure}

The resulting TFSI temperature in the reconstructed domain can be seen in \Cref{fig:ex4.tfsi-temp}. This temperature is positive and larger inside the crack than the surrounding material. We, therefore, expect the crack to open less when the temperature is coupled to the phase-field fracture model 
in addition to the pressure. Furthermore, we note that while the maximum temperature gets smaller with each smaller mesh size, the results are overall comparable and consistent. This difference appears to be driven by the different of domain from to the geometry reconstruction with changing CODs on each mesh.

\begin{figure}
  \centering
  \includegraphics[width=.3\textwidth]{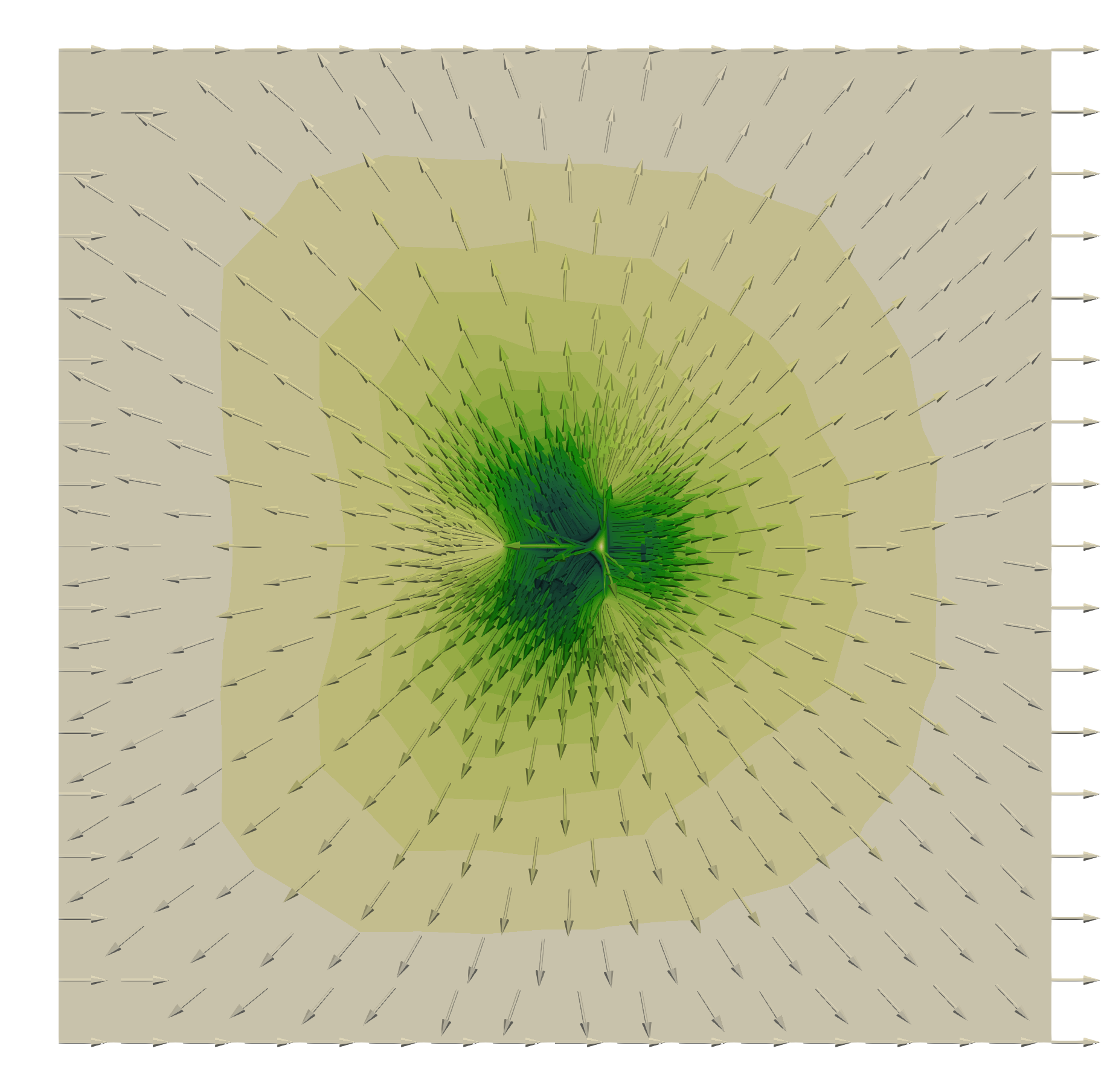}
  \includegraphics[width=.3\textwidth]{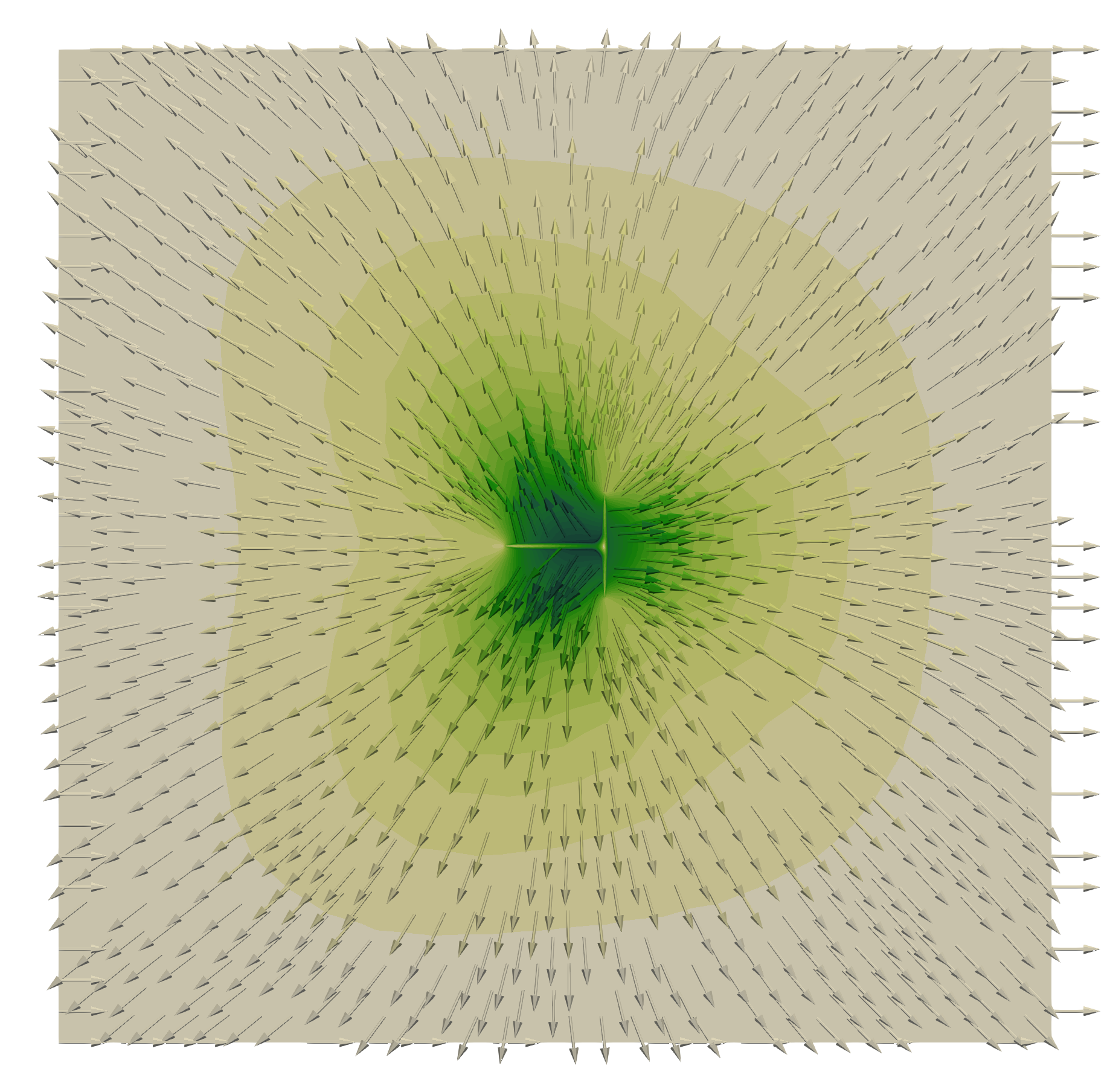}
  \includegraphics[width=.3\textwidth]{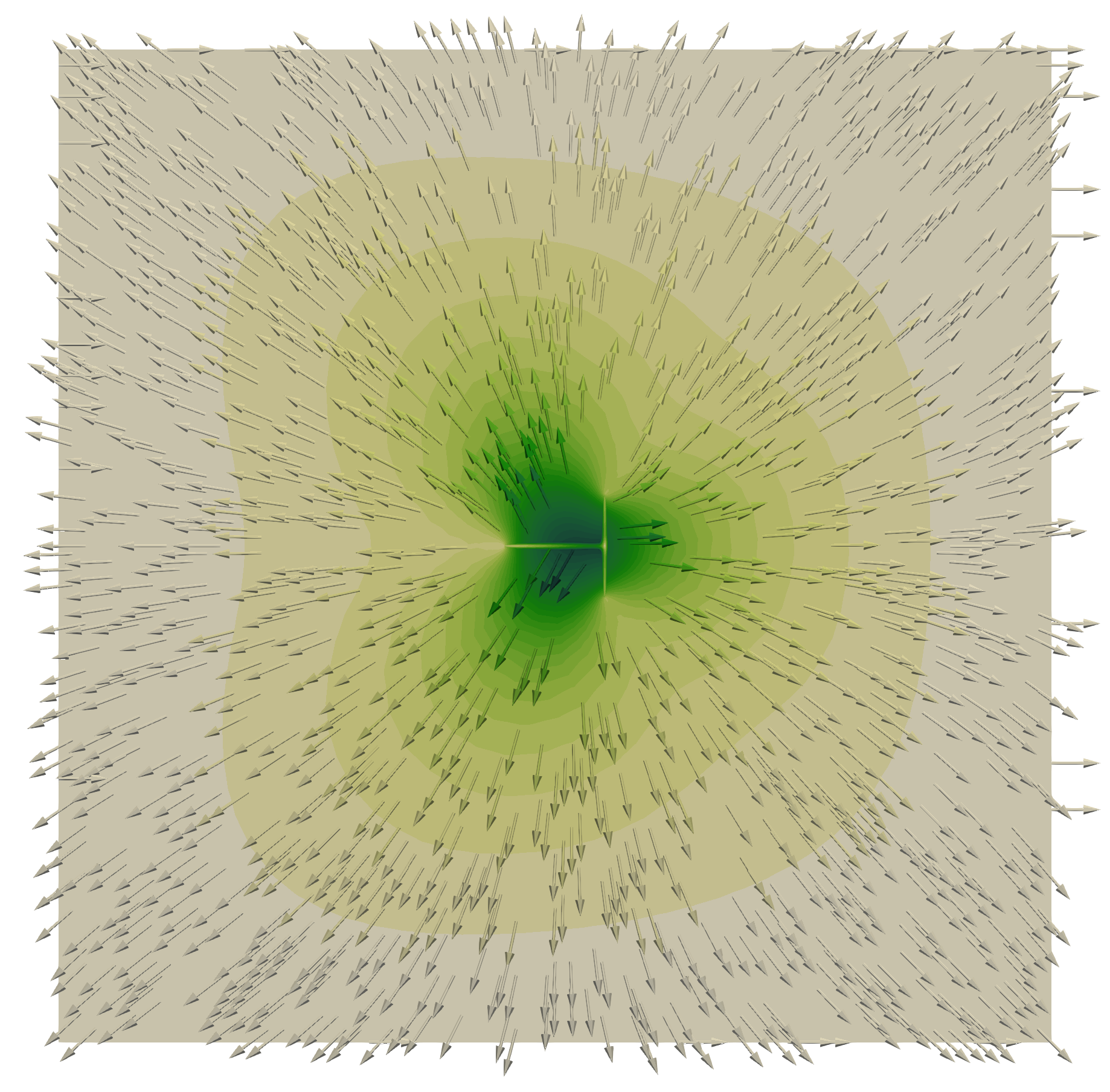}\\[2pt]
  \includegraphics[height=12.5pt]{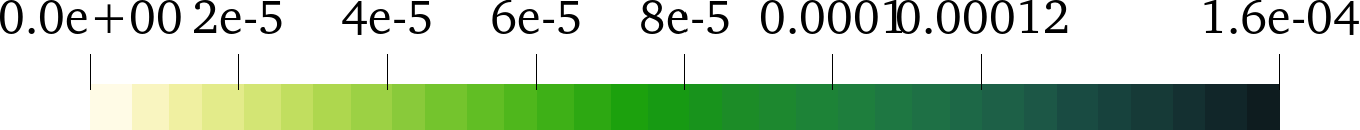}
  \caption{Magnitude and and vector field $\ub_h^{p} - \ub_h^{p,\theta}$ in \hyperref[sec.numerical_approx:subsec.ex4]{Example 4} --- Difference between the phase-field deformation with pressure and combined pressure and temperature coupling. Computed on meshes with $\hmax=0.5, 0.25$ and $0.125$, respectively.}
  \label{fig:ex4.pf-def.diff}
\end{figure}

The difference between the resulting deformations $\ub_h^{p} - \ub_h^{p,\theta}$ is shown in \Cref{fig:ex4.pf-def.diff}. We see that this difference points away from the crack, indicating that the deformation $\ub_h^{p}$ resulting from just the pressure coupling is indeed larger than the deformation $\ub_h^{p,\theta}$ from both pressure and temperature coupling.

\section{Conclusions}
\label{sec.conclusions}

In this work, both a new mathematical model and a new numerical approach for thermo-flow-mechanics-fracture are derived. The key idea is to utilize a phase-field approach for fracture opening and fracture propagation. Having the fracture subdomain at hand, a geometry reconstruction approach is employed. This then results in a mesh with which resolves the boundary between the fluid filled crack and the intact solid domains. This allows us to formulate sharp interface problems for the fracture subdomain and the surrounding medium. The resulting framework is a mixture of interface-capturing and interface-tracking approaches that are conveniently combined. It is substantiated for thermo-flow-mechanics-fracture, which is on the one hand thermo-flow-mechanics (THM) phase-field fracture coupled with thermo fluid-structure interaction (TFSI). The latter is prescribed on moving domains (as the fracture moves) for which the arbitrary Lagrangian-Eulerian (ALE) technique is employed. The governing physics are newly developed, specifically the interface conditions for the phase-field sub-problem. These ingredients yield an overall coupling algorithm with four principle steps after an initialization step: fracture width computation (step 1), re-meshing of the reconstructed subdomains (step 2), solving the TFSI problem (step 3), solving the THM phase-field fracture problem (step 4). The algorithmic details and resulting sub-problems are carefully worked out. 

In order to substantiate our new model and new algorithms, we conducted several numerical experiments. Therein, a key component are mesh refinement studies in which the computational robustness is investigated. This was done for the total crack volume (TCV) and the crack opening displacements (COD) as well as the fracture length. All yielded satisfactory findings in view of the complexity of the problem statement. It should be mentioned that specifically the re-meshing of the crack tips required additional algorithmic developments and obtaining computational convergence was a challenge. From the physics point of view, we emphasized the temperature's influence, where there is agreement in the literature that fractures open due to cold water injection and close due to warm water injection. Single fracture as well as two joining fracture were considered. 

The beauty of our overall framework is that different physics can be easily exchanged as long as the interface conditions are correctly modeled. Therefore, our framework presents the opportunity for future extensions, which could include, for example, two-phase flows, or more complicated mechanics. Furthermore, the extension to three spatial dimensions will remain a challenge due to the geometry reconstruction and requires more extensive future work.

\section*{Data Availability Statement}
The code used to realize the presented results, as well as the resulting raw data, is freely available on github under \url{https://github.com/hvonwah/repro-tfsi-pff} and is archived on zenodo \url{https://doi.org/10.5281/zenodo.13685486}.

\section*{Acknowledgements}
This material is based upon work supported by the National Science Foundation under Grant No.\@{} DMS-1929284 while the HvW was in residence at the Institute for Computational and Experimental Research in Mathematics in Providence, RI, during the Numerical PDEs: Analysis, Algorithms, and Data Challenges program.
SL acknowledges support within his research stay in May 2024 at the Leibniz University Hannover.
The work of S. Lee was partially supported by the U.S.
National Science Foundation Grant DMS-2208402 and by the U.S. Department of Energy, Office of Science,
Energy Earthshots Initiatives under Award Number DE-SC-0024703.


\begin{thebibliography}{10}
\expandafter\ifx\csname url\endcsname\relax
  \def\url#1{\texttt{#1}}\fi
\expandafter\ifx\csname urlprefix\endcsname\relax\def\urlprefix{URL }\fi
\expandafter\ifx\csname href\endcsname\relax
  \def\href#1#2{#2} \def\path#1{#1}\fi

\bibitem{moore2019utah}
J.~Moore, J.~McLennan, R.~Allis, K.~Pankow, S.~Simmons, R.~Podgorney,
  P.~Wannamaker, J.~Bartley, C.~Jones, W.~Rickard, The {Utah} {Frontier}
  {Observatory} for p in {Geothermal} {Energy} ({FORGE}): An international
  laboratory for enhanced geothermal system technology development, in: 44th
  Workshop on Geothermal Reservoir Engineering, Stanford University, 2019, pp.
  11--13.

\bibitem{olasolo2016enhanced}
P.~Olasolo, M.~C. Ju{\'a}rez, M.~P. Morales, I.~A. Liarte, et~al., Enhanced
  geothermal systems ({EGS)}: {A} review, Renewable Sustainable Energy Rev. 56
  (2016) 133--144.
\newblock \href {https://doi.org/10.1016/j.rser.2015.11.031}
  {\path{doi:10.1016/j.rser.2015.11.031}}.

\bibitem{lu2018global}
S.-M. Lu, A global review of enhanced geothermal system ({EGS}), Renewable
  Sustainable Energy Rev. 81 (2018) 2902--2921.
\newblock \href {https://doi.org/10.1016/j.rser.2017.06.097}
  {\path{doi:10.1016/j.rser.2017.06.097}}.

\bibitem{mcclure2014investigation}
M.~W. McClure, R.~N. Horne, An investigation of stimulation mechanisms in
  {Enhanced} {Geothermal} {Systems}, Int. J. Rock Mech. Min. 72 (2014)
  242--260.
\newblock \href {https://doi.org/10.1016/j.ijrmms.2014.07.011}
  {\path{doi:10.1016/j.ijrmms.2014.07.011}}.

\bibitem{caulk2016experimental}
R.~A. Caulk, E.~Ghazanfari, J.~N. Perdrial, N.~Perdrial, Experimental
  investigation of fracture aperture and permeability change within {Enhanced}
  {Geothermal} {Systems}, Geothermics 62 (2016) 12--21.
\newblock \href {https://doi.org/10.1016/j.geothermics.2016.02.003}
  {\path{doi:10.1016/j.geothermics.2016.02.003}}.

\bibitem{zhang2022thermal}
X.~Zhang, Z.~Li, X.~Wang, H.~Wang, B.~Li, Y.~Niu, Thermal effect on the
  fracture behavior of granite using acoustic emission and digital image
  correlation: {An} experimental investigation, Theor. Appl. Fract. Mec. 121
  (2022) 103540.
\newblock \href {https://doi.org/10.1016/j.tafmec.2022.103540}
  {\path{doi:10.1016/j.tafmec.2022.103540}}.

\bibitem{donahue1972crack}
R.~J. Donahue, H.~M. Clark, P.~Atanmo, R.~Kumble, A.~J. McEvily, Crack opening
  displacement and the rate of fatigue crack growth, Int. J. Fract. Mech. 8
  (1972) 209--219.
\newblock \href {https://doi.org/10.1007/BF00703882}
  {\path{doi:10.1007/BF00703882}}.

\bibitem{burdekin1966crack}
F.~M. Burdekin, D.~E.~W. Stone, The crack opening displacement approach to
  fracture mechanics in yielding materials, J. Strain Anal. Eng. Des. 1~(2)
  (1966) 145--153.
\newblock \href {https://doi.org/10.1243/03093247V012145}
  {\path{doi:10.1243/03093247V012145}}.

\bibitem{zhou2022thermal}
L.~Zhou, W.~Gao, L.~Yu, Z.~Zhu, J.~Chen, X.~Wang, Thermal effects on fracture
  toughness of cracked straight-through {Brazilian} disk green sandstone and
  granite, J. Rock Mech. Geotech. Eng. 14~(5) (2022) 1447--1460.
\newblock \href {https://doi.org/10.1016/j.jrmge.2022.02.016}
  {\path{doi:10.1016/j.jrmge.2022.02.016}}.

\bibitem{BourFraMar00}
B.~Bourdin, G.~A. Francfort, J.-J. Marigo, Numerical experiments in revisited
  brittle fracture, J. Mech. Phys. Solids 48~(4) (2000) 797--826.
\newblock \href {https://doi.org/10.1016/S0022-5096(99)00028-9}
  {\path{doi:10.1016/S0022-5096(99)00028-9}}.

\bibitem{KuMue10}
C.~Kuhn, R.~Müller, A continuum phase field model for fracture, Eng. Fract.
  Mech. 77~(18) (2010) 3625--3634, {C}omputational Mechanics in Fracture and
  Damage: A Special Issue in Honor of Prof. Gross.
\newblock \href {https://doi.org/10.1016/j.engfracmech.2010.08.009}
  {\path{doi:10.1016/j.engfracmech.2010.08.009}}.

\bibitem{MieWelHof10a}
C.~Miehe, F.~Welschinger, M.~Hofacker, Thermodynamically consistent phase-field
  models of fracture: {V}ariational principles and multi-field {FE}
  implementations, Internat. J. Numer. Methods Engrg. 83~(10) (2010)
  1273--1311.
\newblock \href {https://doi.org/10.1002/nme.2861}
  {\path{doi:10.1002/nme.2861}}.

\bibitem{MieWelHof10b}
C.~Miehe, M.~Hofacker, F.~Welschinger, A phase field model for rate-independent
  crack propagation: {R}obust algorithmic implementation based on operator
  splits, Comput. Methods Appl. Mech. Engrg. 199 (2010) 2765--2778.
\newblock \href {https://doi.org/10.1016/j.cma.2010.04.011}
  {\path{doi:10.1016/j.cma.2010.04.011}}.

\bibitem{BoVeScoHuLa12}
M.~J. Borden, C.~V. Verhoosel, M.~A. Scott, T.~J.~R. Hughes, C.~M. Landis, A
  phase-field description of dynamic brittle fracture, Comput. Methods Appl.
  Mech. Engrg. 217 (2012) 77--95.
\newblock \href {https://doi.org/10.1016/j.cma.2012.01.008}
  {\path{doi:10.1016/j.cma.2012.01.008}}.

\bibitem{AmGeraLoren15}
M.~Ambati, T.~Gerasimov, L.~De~Lorenzis, A review on phase-field models of
  brittle fracture and a new fast hybrid formulation, Comput. Mech. 55~(2)
  (2015) 383--405.
\newblock \href {https://doi.org/10.1007/s00466-014-1109-y}
  {\path{doi:10.1007/s00466-014-1109-y}}.

\bibitem{ARRIAGA201833}
M.~Arriaga, H.~Waisman, Stability analysis of the phase-field method for
  fracture with a general degradation function and plasticity induced crack
  generation, Mech. Mater. 116 (2018) 33--48, iUTAM Symposium on Dynamic
  Instabilities in Solids.
\newblock \href {https://doi.org/10.1016/j.mechmat.2017.04.003}
  {\path{doi:10.1016/j.mechmat.2017.04.003}}.

\bibitem{SARGADO2018458}
J.~M. Sargado, E.~Keilegavlen, I.~Berre, J.~M. Nordbotten, High-accuracy
  phase-field models for brittle fracture based on a new family of degradation
  functions, J. Mech. Phys. Solids 111 (2018) 458--489.
\newblock \href {https://doi.org/10.1016/j.jmps.2017.10.015}
  {\path{doi:10.1016/j.jmps.2017.10.015}}.

\bibitem{WheWiLee20}
M.~F. Wheeler, T.~Wick, S.~Lee, {IPACS: Integrated Phase-Field Advanced Crack
  Propagation Simulator. {An} adaptive, parallel, physics-based-discretization
  phase-field framework for fracture propagation in porous media}, Comput.
  Methods Appl. Mech. Engrg. 367 (2020) 113124.
\newblock \href {https://doi.org/10.1016/j.cma.2020.113124}
  {\path{doi:10.1016/j.cma.2020.113124}}.

\bibitem{BourFraMar08}
B.~Bourdin, G.~A. Francfort, J.-J. Marigo, The variational approach to
  fracture, J. Elasticity 91~(1--3) (2008) 1--148.
\newblock \href {https://doi.org/10.1007/s10659-007-9107-3}
  {\path{doi:10.1007/s10659-007-9107-3}}.

\bibitem{WNN20}
J.-Y. Wu, V.~P. Nguyen, C.~T. Nguyen, D.~Sutula, S.~Sinaie, S.~P.~A. Bordas,
  Phase-field modeling of fracture, Advances in Applied Mechanics, Elsevier,
  2020, Ch.~1, pp. 1--183.
\newblock \href {https://doi.org/10.1016/bs.aams.2019.08.001}
  {\path{doi:10.1016/bs.aams.2019.08.001}}.

\bibitem{Wi20}
T.~Wick, Multiphysics Phase-Field Fracture, Vol.~28 of Radon Series on
  Computational and Applied Mathematics, De Gruyter, Berlin, Boston, 2020.
\newblock \href {https://doi.org/10.1515/9783110497397}
  {\path{doi:10.1515/9783110497397}}.

\bibitem{HEIDER2021107881}
Y.~Heider, A review on phase-field modeling of hydraulic fracturing, Eng.
  Fract. Mech. 253 (2021) 107881.
\newblock \href {https://doi.org/10.1016/j.engfracmech.2021.107881}
  {\path{doi:10.1016/j.engfracmech.2021.107881}}.

\bibitem{DiLiWiTy22}
P.~Diehl, R.~Lipton, T.~Wick, M.~Tyagi, A comparative review of peridynamics
  and phase-field models for engineering fracture mechanics, Comput. Mech. 69
  (2022) 1259--1293.
\newblock \href {https://doi.org/10.1007/s00466-022-02147-0}
  {\path{doi:10.1007/s00466-022-02147-0}}.

\bibitem{LeeWheWi16}
S.~Lee, M.~F. Wheeler, T.~Wick, Pressure and fluid-driven fracture propagation
  in porous media using an adaptive finite element phase field model, Comput.
  Methods Appl. Mech. Engrg. 305 (2016) 111--132.
\newblock \href {https://doi.org/10.1016/j.cma.2016.02.037}
  {\path{doi:10.1016/j.cma.2016.02.037}}.

\bibitem{lee2016phase}
S.~Lee, A.~Mikeli{\'c}, M.~F. Wheeler, T.~Wick, Phase-field modeling of
  proppant-filled fractures in a poroelastic medium, Comput. Methods Appl.
  Mech. Engrg. 312 (2016) 509--541.
\newblock \href {https://doi.org/10.1016/j.cma.2016.02.008}
  {\path{doi:10.1016/j.cma.2016.02.008}}.

\bibitem{MDB99}
N.~Moës, J.~Dolbow, T.~Belytschko, A finite element method for crack growth
  without remeshing, Internat. J. Numer. Methods Engrg. 46~(1) (1999) 131--150.
\newblock \href
  {https://doi.org/10.1002/(sici)1097-0207(19990910)46:1<131::aid-nme726>3.0.co;2-j}
  {\path{doi:10.1002/(sici)1097-0207(19990910)46:1<131::aid-nme726>3.0.co;2-j}}.

\bibitem{zi2003new}
G.~Zi, T.~Belytschko, New crack-tip elements for {XFEM} and applications to
  cohesive cracks, Internat. J. Numer. Methods Engrg. 57~(15) (2003)
  2221--2240.
\newblock \href {https://doi.org/10.1002/nme.849} {\path{doi:10.1002/nme.849}}.

\bibitem{duarte2001generalized}
C.~A. Duarte, O.~N. Hamzeh, T.~J. Liszka, W.~W. Tworzydlo, A generalized finite
  element method for the simulation of three-dimensional dynamic crack
  propagation, Comput. Methods Appl. Mech. Engrg. 190~(15-17) (2001)
  2227--2262.
\newblock \href {https://doi.org/10.1016/S0045-7825(00)00233-4}
  {\path{doi:10.1016/S0045-7825(00)00233-4}}.

\bibitem{FB10}
T.-P. Fries, T.~Belytschko, The extended/generalized finite element method:
  {An} overview of the method and its applications, Internat. J. Numer. Methods
  Engrg. 84~(3) (2010) 253--304.
\newblock \href {https://doi.org/10.1002/nme.2914}
  {\path{doi:10.1002/nme.2914}}.

\bibitem{LWW17}
S.~Lee, M.~F. Wheeler, T.~Wick, Iterative coupling of flow, geomechanics and
  adaptive phase-field fracture including level-set crack width approaches, J.
  Comput. Appl. Math. 314 (2017) 40--60.
\newblock \href {https://doi.org/10.1016/j.cam.2016.10.022}
  {\path{doi:10.1016/j.cam.2016.10.022}}.

\bibitem{YNK20}
K.~Yoshioka, D.~Naumov, O.~Kolditz, On crack opening computation in variational
  phase-field models for fracture, Comput. Methods Appl. Mech. Engrg. 369
  (2020) 113210.
\newblock \href {https://doi.org/10.1016/j.cma.2020.113210}
  {\path{doi:10.1016/j.cma.2020.113210}}.

\bibitem{WaWi24_RINAM}
H.~von Wahl, T.~Wick, A coupled high-accuracy phase-field fluid–structure
  interaction framework for stokes fluid-filled fracture surrounded by an
  elastic medium, Results Appl. Math. 22 (2024) 100455.
\newblock \href {https://doi.org/10.1016/j.rinam.2024.100455}
  {\path{doi:10.1016/j.rinam.2024.100455}}.

\bibitem{WaWi23_CMAME}
H.~von Wahl, T.~Wick, A high-accuracy framework for phase-field fracture
  interface reconstructions with application to {Stokes} fluid-filled fracture
  surrounded by an elastic medium, Comput. Methods Appl. Mech. Engrg. 415
  (2023) 116202.
\newblock \href {https://doi.org/10.1016/j.cma.2023.116202}
  {\path{doi:10.1016/j.cma.2023.116202}}.

\bibitem{NoiiWi19}
N.~Noii, T.~Wick, A phase-field description for pressurized and non-isothermal
  propagating fractures, Comput. Methods Appl. Mech. Engrg. 351 (2019)
  860--890.
\newblock \href {https://doi.org/10.1016/j.cma.2019.03.058}
  {\path{doi:10.1016/j.cma.2019.03.058}}.

\bibitem{HEIDER2018116}
Y.~Heider, S.~Reiche, P.~Siebert, B.~Markert, Modeling of hydraulic fracturing
  using a porous-media phase-field approach with reference to experimental
  data, Eng. Fract. Mech. 202 (2018) 116--134.
\newblock \href {https://doi.org/10.1016/j.engfracmech.2018.09.010}
  {\path{doi:10.1016/j.engfracmech.2018.09.010}}.

\bibitem{NgHeiMa23}
C.-L. Nguyen, Y.~Heider, B.~Markert, A non-isothermal phase-field hydraulic
  fracture modeling in saturated porous media with convection-dominated heat
  transport, Acta Geotech. 50~(6) (2023) 821--833.
\newblock \href {https://doi.org/10.1007/s11440-023-01905-5}
  {\path{doi:10.1007/s11440-023-01905-5}}.

\bibitem{SUH2021114182}
H.~S. Suh, W.~Sun, Asynchronous phase field fracture model for porous media
  with thermally non-equilibrated constituents, Comput. Methods Appl. Mech.
  Engrg. 387 (2021) 114182.
\newblock \href {https://doi.org/10.1016/j.cma.2021.114182}
  {\path{doi:10.1016/j.cma.2021.114182}}.

\bibitem{Dai2024}
Y.~Dai, B.~Hou, S.~Lee, T.~Wick, A thermal--hydraulic--mechanical--chemical
  coupling model for acid fracture propagation based on a phase-field method,
  Rock Mech. Rock Eng. (2024).
\newblock \href {https://doi.org/10.1007/s00603-024-03769-x}
  {\path{doi:10.1007/s00603-024-03769-x}}.

\bibitem{LIU2024117165}
Y.~Liu, K.~Yoshioka, T.~You, H.~Li, F.~Zhang, A phase-field fracture model in
  thermo-poro-elastic media with micromechanical strain energy degradation,
  Comput. Methods Appl. Mech. Engrg. 429 (2024) 117165.
\newblock \href {https://doi.org/10.1016/j.cma.2024.117165}
  {\path{doi:10.1016/j.cma.2024.117165}}.

\bibitem{lee4920845phase}
S.~Lee, M.~Wheeler, T.~Wick, \href{https://ssrn.com/abstract=4920845}{A
  phase-field diffraction model for thermo-hydro-mechanical propagating
  fractures} (Aug. 2024).
\newline\urlprefix\url{https://ssrn.com/abstract=4920845}

\bibitem{LoBo96}
S.~A. Lorca, J.~L. Boldrini, Stationary solutions for generalized {Boussinesq}
  models, J. Differ. Equ. 124~(2) (1996) 389--406.
\newblock \href {https://doi.org/10.1006/jdeq.1996.0016}
  {\path{doi:10.1006/jdeq.1996.0016}}.

\bibitem{FARHAT1991349}
C.~Farhat, K.~C. Park, Y.~Dubois-Pelerin, An unconditionally stable staggered
  algorithm for transient finite element analysis of coupled thermoelastic
  problems, Comput. Methods Appl. Mech. Engrg. 85~(3) (1991) 349--365.
\newblock \href {https://doi.org/10.1016/0045-7825(91)90102-C}
  {\path{doi:10.1016/0045-7825(91)90102-C}}.

\bibitem{Coussy2004}
O.~Coussy, Poromechanics, Wiley, 2004.
\newblock \href {https://doi.org/10.1002/0470092718}
  {\path{doi:10.1002/0470092718}}.

\bibitem{mayeli2021buoyancy}
P.~Mayeli, G.~J. Sheard, Buoyancy-driven flows beyond the {Boussinesq}
  approximation: {A} brief review, Int. Commun. Heat Mass 125 (2021) 105316.
\newblock \href {https://doi.org/10.1016/j.icheatmasstransfer.2021.105316}
  {\path{doi:10.1016/j.icheatmasstransfer.2021.105316}}.

\bibitem{MiWheWi13a}
A.~Mikeli\'c, M.~F. Wheeler, T.~Wick,
  \href{www.oden.utexas.edu/media/reports/2013/1315.pdf}{A phase-field approach
  to the fluid filled fracture surrounded by a poroelastic medium}, iCES Report
  13-15 (Jun. 2013).
\newline\urlprefix\url{www.oden.utexas.edu/media/reports/2013/1315.pdf}

\bibitem{MWW19}
A.~Mikeli{\'{c}}, M.~F. Wheeler, T.~Wick, Phase-field modeling through
  iterative splitting of hydraulic fractures in a poroelastic medium, GEM -
  Int. J. Geomath. 10~(1) (Jan. 2019).
\newblock \href {https://doi.org/10.1007/s13137-019-0113-y}
  {\path{doi:10.1007/s13137-019-0113-y}}.

\bibitem{AmTo90}
L.~Ambrosio, V.~M. Tortorelli, Approximation of functionals depending on jumps
  by elliptic functionals via $\gamma$-convergence, Comm. Pure Appl. Math.
  43~(8) (1990) 999--1036.
\newblock \href {https://doi.org/10.1002/cpa.3160430805}
  {\path{doi:10.1002/cpa.3160430805}}.

\bibitem{AmTo92}
L.~Ambrosio, V.~M. Tortorelli, On the approximation of free discontinuity
  problems, Boll. Un. Mat. Ital. 6 (1992) 105--123.

\bibitem{HWW15}
T.~Heister, M.~F. Wheeler, T.~Wick, A primal-dual active set method and
  predictor-corrector mesh adaptivity for computing fracture propagation using
  a phase-field approach, Comput. Methods Appl. Mech. Engrg. 290 (2015)
  466--495.
\newblock \href {https://doi.org/10.1016/j.cma.2015.03.009}
  {\path{doi:10.1016/j.cma.2015.03.009}}.

\bibitem{Wic17}
T.~Wick, An error-oriented {Newton}/inexact augmented {Lagrangian} approach for
  fully monolithic phase-field fracture propagation, SIAM J. Sci. Comput.
  39~(4) (2017) B589--B617.
\newblock \href {https://doi.org/10.1137/16m1063873}
  {\path{doi:10.1137/16m1063873}}.

\bibitem{KMW23}
L.~Kolditz, K.~Mang, T.~Wick, A modified combined active-set {Newton} method
  for solving phase-field fracture into the monolithic limit, Comput. Methods
  Appl. Mech. Engrg. 414 (2023) 116170.
\newblock \href {https://doi.org/10.1016/j.cma.2023.116170}
  {\path{doi:10.1016/j.cma.2023.116170}}.

\bibitem{MWT15}
A.~Mikeli{\'{c}}, M.~F. Wheeler, T.~Wick, A quasi-static phase-field approach
  to pressurized fractures, Nonlinearity 28~(5) (2015) 1371--1399.
\newblock \href {https://doi.org/10.1088/0951-7715/28/5/1371}
  {\path{doi:10.1088/0951-7715/28/5/1371}}.

\bibitem{TrSeNg13}
D.~Tran, A.~T. Settari, L.~Nghiem, Predicting growth and decay of
  hydraulic-fracture width in porous media subjected to isothermal and
  nonisothermal flow, SPE J. 18~(4) (2013) 781--794.
\newblock \href {https://doi.org/10.2118/162651-PA}
  {\path{doi:10.2118/162651-PA}}.

\bibitem{CHUKWUDOZIE2019957}
C.~Chukwudozie, B.~Bourdin, K.~Yoshioka, A variational phase-field model for
  hydraulic fracturing in porous media, Comput. Methods Appl. Mech. Engrg. 347
  (2019) 957--982.
\newblock \href {https://doi.org/10.1016/j.cma.2018.12.037}
  {\path{doi:10.1016/j.cma.2018.12.037}}.

\bibitem{HrTu06a}
J.~Hron, S.~Turek, A monolithic {FEM}/Multigrid solver for {ALE} formulation of
  fluid structure with application in biomechanics, Vol.~53, Springer, Berlin,
  Heidelberg, 2006, pp. 146--170.
\newblock \href {https://doi.org/10.1007/3-540-34596-5_7}
  {\path{doi:10.1007/3-540-34596-5_7}}.

\bibitem{Du07}
T.~Dunne, Adaptive finite element approximation of fluid-structure interaction
  based on {E}ulerian and arbitrary {L}agrangian-{E}ulerian variational
  formulations, Ph.D. thesis, University of Heidelberg (2007).
\newblock \href {https://doi.org/10.11588/heidok.00007944}
  {\path{doi:10.11588/heidok.00007944}}.

\bibitem{Wi11_phd}
T.~Wick, Adaptive {F}inite {E}lement {S}imulation of {F}luid-{S}tructure
  {I}nteraction with {A}pplication to {H}eart-{V}alve {D}ynamics, Ph.D. thesis,
  University of Heidelberg (2011).
\newblock \href {https://doi.org/10.11588/heidok.00012992}
  {\path{doi:10.11588/heidok.00012992}}.

\bibitem{Ri17_fsi}
T.~Richter, Fluid-structure interactions: {Models}, analysis, and finite
  elements, Springer, Cham, 2017.
\newblock \href {https://doi.org/10.1007/978-3-319-63970-3}
  {\path{doi:10.1007/978-3-319-63970-3}}.

\bibitem{HuLiZi81}
T.~J.~R. Hughes, W.~K. Liu, T.~Zimmermann, {L}agrangian-{E}ulerian finite
  element formulation for incompressible viscous flows, Comput. Methods Appl.
  Mech. Engrg. 29 (1981) 329--349.
\newblock \href {https://doi.org/10.1016/0045-7825(81)90049-9}
  {\path{doi:10.1016/0045-7825(81)90049-9}}.

\bibitem{DoGiuHa82}
J.~Donea, S.~Giuliani, J.~P. Halleux, An arbitrary {Lagrangian}-{Eulerian}
  finite element method for transient dynamic fluid-structure interactions,
  Comput. Methods Appl. Mech. Engrg. 33 (1982) 689--723.
\newblock \href {https://doi.org/10.1016/0045-7825(82)90128-1}
  {\path{doi:10.1016/0045-7825(82)90128-1}}.

\bibitem{RiWi10}
T.~Richter, T.~Wick, Finite elements for fluid-structure interaction in {ALE}
  and fully {E}ulerian coordinates, Comput. Methods Appl. Mech. Engrg. 199
  (2010) 2633--2642.
\newblock \href {https://doi.org/10.1016/j.cma.2010.04.016}
  {\path{doi:10.1016/j.cma.2010.04.016}}.

\bibitem{Sch97}
J.~Sch\"oberl, {NETGEN} an advancing front {2D}/{3D}-mesh generator based on
  abstract rules, Comput. Vis. Sci. 1~(1) (1997) 41--52.
\newblock \href {https://doi.org/10.1007/s007910050004}
  {\path{doi:10.1007/s007910050004}}.

\bibitem{Sch14}
J.~Sch\"oberl, C++11 implementation of finite elements in {}{NGSolve}{}, Tech.
  Rep. ASC Report No. 30/2014 (Sep. 2014).

\bibitem{LHPvW21}
C.~Lehrenfeld, F.~Heimann, J.~Preu\ss, H.~von Wahl,
  \href{github.com/ngsxfem/ngsxfem}{\texttt{ngsxfem}: Add-on to ngsolve for
  geometrically unfitted finite element discretizations}, J. Open Source Softw.
  6~(64) (2021) 3237.
\newblock \href {https://doi.org/10.21105/joss.03237}
  {\path{doi:10.21105/joss.03237}}.
\newline\urlprefix\url{github.com/ngsxfem/ngsxfem}

\bibitem{Sne46}
I.~N. Sneddon, The distribution of stress in the neighbourhood of a crack in an
  elastic solid, Proc. R. Soc. A 187~(1009) (1946) 229--260.
\newblock \href {https://doi.org/10.1098/rspa.1946.0077}
  {\path{doi:10.1098/rspa.1946.0077}}.

\bibitem{SneddLow69}
I.~N. Sneddon, M.~Lowengrub, Crack problems in the classical theory of
  elasticity, SIAM series in Applied Mathematics, John Wiley and Sons,
  Philadelphia, 1969.

\bibitem{morrow2001permeability}
C.~A. Morrow, D.~E. Moore, D.~A. Lockner, Permeability reduction in granite
  under hydrothermal conditions, J. Geophys. Res. Solid Earth 106~(B12) (2001)
  30551--30560.
\newblock \href {https://doi.org/10.1029/2000JB000010}
  {\path{doi:10.1029/2000JB000010}}.

\bibitem{yasuhara2006evolution}
H.~Yasuhara, A.~Polak, Y.~Mitani, A.~S. Grader, P.~M. Halleck, D.~Elsworth,
  Evolution of fracture permeability through fluid--rock reaction under
  hydrothermal conditions, Earth Planet. Sc. Lett. 244~(1-2) (2006) 186--200.
\newblock \href {https://doi.org/10.1016/j.epsl.2006.01.046}
  {\path{doi:10.1016/j.epsl.2006.01.046}}.

\bibitem{hardin1982measuring}
E.~Hardin, N.~Barton, M.~Voegele, M.~Board, R.~Lingle, H.~Pratt, W.~Ubbes,
  Measuring the thermomechanical and transport properties of a rockmass using
  the heated block test, in: ARMA US Rock Mechanics/Geomechanics Symposium,
  ARMA, 1982, pp. ARMA--82.

\bibitem{rutqvist2008analysis}
J.~Rutqvist, B.~Freifeld, K.-B. Min, D.~Elsworth, Y.~Tsang, Analysis of
  thermally induced changes in fractured rock permeability during 8 years of
  heating and cooling at the yucca mountain drift scale test, Int. J. Rock
  Mech. Min. 45~(8) (2008) 1373--1389.
\newblock \href {https://doi.org/10.1016/j.ijrmms.2008.01.016}
  {\path{doi:10.1016/j.ijrmms.2008.01.016}}.


\end{thebibliography}
\end{document}